\documentclass{amsart}
\usepackage{amsfonts}
\usepackage{fancyhdr}
\usepackage[latin1]{inputenc}
\usepackage{amsmath}
\usepackage{thmdefs}
\usepackage{natbib}
\usepackage{graphicx}
\usepackage{hyperref}
\usepackage{amssymb}
\usepackage{dsfont}

\theoremstyle{definition} \theoremstyle{remark}
\numberwithin{equation}{section}

\renewcommand{\cite}{\citet}

\begin{document}

\author{Michel BRONIATOWSKI$^{*}$ and Amor KEZIOU$^{**}$}
\address{$^{*}$LSTA-Universit\'{e} Paris 6. e-mail: michel.broniatowski@upmc.fr
\newline
$^{**}$Laboratoire de Math\'ematiques,
Universit\'e de Reims  and LSTA-Universit\'e Paris 6. \\
 e-mail: amor.keziou@upmc.fr}
\date{March 2010}
\title[Minimization of divergences on sets of signed measures]{Minimization of divergences on sets of signed measures
\footnote{\textbf{This is an an electronic reprint of the original
article published in Studia Sci. Math. Hungar., 2006, Vol. 43, No.
4, 403--442. This reprint differs from the original in pagination
and typographic detail.}}}

\maketitle

\begin{abstract}
We consider the minimization problem of $\phi$-divergences between
a given probability measure $P$ and subsets $\Omega$ of the vector
space $\mathcal{M}_\mathcal{F}$ of all signed finite measures
which integrate a given class $\mathcal{F}$ of bounded or
unbounded measurable functions. The vector space
$\mathcal{M}_\mathcal{F}$ is endowed with the weak topology
induced by the class $\mathcal{F}\cup \mathcal{B}_b$ where
$\mathcal{B}_b$ is the class of all bounded measurable functions.
We treat the problems of existence and characterization of the
$\phi$-projections of $P$ on $\Omega$. We consider also the dual
equality and the dual attainment problems when $\Omega$ is defined
by linear constraints.
 \vspace{2mm} \\ Key words: Minimum Divergences; Maximum Entropy; Convex
 Programming;
Moment Problem; Empirical Likelihood; Convex Distances;
Fenchel Duality.\\
\end{abstract}

 \subjclass{MSC (2000) Classification:
49A55; 49A40; 46A20; 46A05; 62E20.}

\tableofcontents


\section{Introduction and notation}
\noindent Let $\left(\mathcal{X},\mathcal{B}\right)$ be a
measurable space and $P$ be a given reference probability measure
(p.m.) on $\left(\mathcal{X},\mathcal{B}\right)$. Denote
$\mathcal{M}$ the real vector space of all signed finite measures
on $\left(\mathcal{X},\mathcal{B}\right)$ and $\mathcal{M}(P)$ the
vector subspace of all signed finite measures absolutely
continuous (a.c) with respect to (w.r.t.) $P$. Denote also
$\mathcal{M}^1$ the set of all p.m.'s on
$\left(\mathcal{X},\mathcal{B}\right)$ and $\mathcal{M}^1(P)$ the
subset of all p.m.'s a.c w.r.t. $P$. Let $\varphi$ be a
proper\footnote{We say a function is proper if its domain is non
void.} closed\footnote{The closedness of $\varphi$ means that if
$a_\varphi$ or $b_\varphi$ are finite numbers then $\varphi(x)$
tends to $\varphi(a_\varphi)$ or $\varphi(b_\varphi)$ when $x
\downarrow a_\varphi$ or $x \uparrow b_\varphi$, respectively.}
convex function from $]-\infty, +\infty[$ to $[0,+\infty]$ with
$\varphi(1)=0$ and such that its domain $\text{dom}\varphi :=
\left\{ x\in\mathbb{R} \text{ such that } \varphi (x)< \infty
\right\}$ is an interval with endpoints $a_\varphi < 1 <
b_\varphi$ (which may be finite or infinite). For any signed
finite measure $Q$ in $\mathcal{M}(P)$, the $\phi$-divergence
between   $Q$ and  $P$ is defined by
\begin{equation}\label{divRusch}
\phi(Q,P):=\int_\mathcal{X}
\varphi\left(\frac{dQ}{dP}(x)\right)~dP(x).
\end{equation}
When $Q$ is not a.c. w.r.t. $P$, we set $\phi(Q,P)=+\infty$. The
$\phi$-divergences between p.m.'s were  introduced by
\cite{Csiszar1963} as ``$f$-divergences''. The definition of
$\phi$-divergences of
 \cite{Csiszar1963} between p.m.'s   requires a common dominating
$\sigma$-finite measure, noted $\lambda$, for $Q$ and $P$. Note
that the two definitions of $\phi-$divergences coincide on the set
of all p.m.'s a.c w.r.t. $P$ and dominated by $\lambda$. The
$\phi$-divergences between any signed finite measure $Q$ and a
p.m. $P$ were introduced by \cite{CsiszarGamboaGassiat1999}; they
gave the following definition
\begin{equation}\label{def div CGC}
\phi(Q,P):=\int
\varphi(q)~dP+b\sigma^+_Q(\mathcal{X})-a\sigma^-_Q(\mathcal{X}),
\end{equation}
where $a:=\lim_{x\to -\infty}\varphi(x)/x$, $b:=\lim_{x\to
+\infty}\varphi(x)/x$ and
\begin{equation*}
Q= qP+\sigma_Q, \quad \sigma_Q=\sigma_Q^+-\sigma_Q^-
\end{equation*}
is the Lebesgue decomposition of $Q$, and the Jordan decomposition
of the singular part  $\sigma_Q$, respectively. The definitions
(\ref{divRusch}) and (\ref{def div CGC}) coincide when $Q$ is a.c.
w.r.t. $P$ or when $a=-\infty$ or $b=+\infty$. Since we will
consider optimization of $Q\mapsto \phi(Q,P)$ on  sets of signed
finite measures a.c. w.r.t. $P$, it is more adequate
for our sake to use the definition (\ref{divRusch}).\\

\noindent For all p.m. $P$, the mappings $Q\in
\mathcal{M}\mapsto\phi(Q,P)$ are convex and take nonnegative
values. When $Q=P$ then $\phi(Q,P)=0$. Furthermore, if the
function $x\mapsto\varphi(x)$ is strictly convex on a neighborhood
of $x=1$, then the following basic property holds
\begin{equation}
\phi(Q,P)=0~\text{ if and only if }~Q=P.  \label{p.f.}
\end{equation}
All these properties are presented in \cite{Csiszar1963},
\cite{Csiszar1967b},  \cite{Csiszar1967} and
\cite{Liese-Vajda1987} chapter 1, for $\phi$-divergences defined
on the set of all p.m.'s $\mathcal{M}^1$. When the
$\phi$-divergences are defined on $\mathcal{M}$, then the same
properties hold. \\
\\
When defined on $\mathcal{M}^1$, the Kullback-Leibler $(KL)$,
modified Kullback-Leibler $(KL_m)$, $\chi^{2}$, modified
$\chi^{2}$ $(\chi_{m}^{2})$, Hellinger $(H)$, and $L_{1}$
divergences are respectively associated to the convex functions
$\varphi(x)=x\log x-x+1$, $\varphi(x)=-\log x+x-1$,
$\varphi(x)=\frac{1}{2}{(x-1)}^{2}$,
$\varphi(x)=\frac{1}{2}{(x-1)}^{2}/x$,
$\varphi(x)=2{(\sqrt{x}-1)}^{2}$ and $\varphi(x)=\left\vert
x-1\right\vert$. All those divergences except the $L_{1}$ one,
belong to the class of power divergences introduced in
\cite{Cressie-Read1984} (see also \cite{Liese-Vajda1987} chapter
2). They are defined through the class of convex functions
\begin{equation}\label{gamma convex functions}
x\in ]0,+\infty[
\mapsto\varphi_{\gamma}(x):=\frac{x^{\gamma}-\gamma
x+\gamma-1}{\gamma(\gamma-1)}
\end{equation}
if $\gamma\in\mathbb{R}\setminus \left\{0,1\right\}$,
$\varphi_{0}(x):=-\log x+x-1$ and $\varphi_{1}(x):=x\log x-x+1$.
(For all $\gamma\in\mathbb{R}$, we define
$\varphi_\gamma(0):=\lim_{x\downarrow 0}\varphi_\gamma (x)$). So,
the $KL-$divergence is associated to $\varphi_1$, the $KL_m$ to
$\varphi_0$, the $\chi^2$ to $\varphi_2$, the $\chi^2_m$ to
$\varphi_{-1}$ and the Hellinger distance to $\varphi_{1/2}$.\\
\\
\noindent The Kullback-Leibler divergence ($KL$-divergence) is
sometimes called Boltzmann Shannon relative entropy. It appears in
the domain of large deviations and it is frequently used  for
reconstruction of laws, and in particular in the classical moment
problem (see e.g. \cite{CsiszarGamboaGassiat1999} and the
references therein). The modified Kullback-Leibler divergence
($KL_m$-divergence) is sometimes called Burg relative entropy. It
is frequently used in Statistics and it leads to efficient methods
in  statistical estimation and tests problems; in fact, the
celebrate ``maximum likelihood'' method can be seen as an
optimization problem of the $KL_m$-divergence between the discrete
or continuous parametric model and the empirical measure
associated to the data; see \cite{Keziou2003} and
\cite{Broniatowski-Keziou2003}. On the other hand, the recent
``empirical likelihood'' method can also be seen as an
optimization problem of the $KL_m$-divergence between some set of
measures satisfying some linear constraints and the empirical
measure associated to the data; see \cite{Owen2001} and the
references therein, \cite{Bertail2003a}, \cite{Bertail2003}
 and \cite{Broniatowski-Keziou2004}. The Hellinger divergence is also
used in Statistics, it leads to robust statistical methods in
parametric and semi-parametric models; see \cite{Beran1977},
\cite{Lindsay1994}, \cite{JimenezShao2001} and
\cite{Broniatowski-Keziou2004}.\\
\\
\noindent We extend the definition of the power divergences
functions $Q\in \mathcal{M}^1\mapsto \phi_{\gamma}(Q,P)$ onto the
whole vector space of signed finite measures $\mathcal{M}$ via the
extension of the definition of the convex functions
$\varphi_\gamma$ : For all $\gamma \in\mathbb{R}$ such that the
function $x\mapsto \varphi_{\gamma }(x)$ is not defined on
$]-\infty, 0[$ or defined but not convex on whole $\mathbb{R}$, we
extend its definition as follows
\begin{equation} \label{gamma convex functions sur R}
x\in ]-\infty,+\infty[\mapsto\left\{
\begin{array}{lll}
\varphi_\gamma (x) & \text{ if
}& x\in [0,+\infty[,\\
+\infty & \text{ if } & x\in ]-\infty,0[.
\end{array}\right.
\end{equation}
Note that for the $\chi^2$-divergence for instance,
$\varphi_2(x):=\frac{1}{2}(x-1)^2$ is
defined and convex on whole $\mathbb{R}$.\\
\\
The conjugate (or Fenchel-Legendre transform) of $\varphi$ will be
denoted $\varphi^*$, i.e.,
\begin{equation}
 t\in \mathbb{R}\mapsto \varphi^*
 (t):=\sup_{x\in\mathbb{R}}\left\{tx-\varphi(x)\right\},
\end{equation}
and the endpoints of $\text{dom}\varphi^*$ (the domain of
$\varphi^*$) will be denoted $a_{\varphi^*}$ and $b_{\varphi^*}$
with $a_{\varphi^*}\leq b_{\varphi^*}$. Note that $\varphi^*$ is
proper closed convex function. In particular, $a_{\varphi^*} < 0 <
b_{\varphi^*}$, $\varphi^*(0)=0$ and
\begin{equation}\label{domain de varphi* }
 a_{\varphi^*}=\lim_{y\to -\infty}\frac{\varphi(y)}{y}, \quad b_{\varphi^*}=
 \lim_{y\to +\infty}\frac{\varphi(y)}{y}.
\end{equation}
By the closedness of $\varphi$, the conjugate $\varphi^{**}$ of
$\varphi^*$ coincides with $\varphi$, i.e.,
\begin{equation}
\varphi^{**}(t):=
\sup_{x\in\mathbb{R}}\left\{tx-\varphi^*(x)\right\}=\varphi(t),
~\text{ for all } t\in\mathbb{R}.
\end{equation}
For the proper convex functions defined on $\mathbb{R}$ (endowed
with the usual topology), the lower semi-continuity\footnote{We
say a function $\varphi$ is lower semi-continuous if the level
sets $\left\{x\text{ such that } \varphi(x)\leq \alpha\right\}$,
$\alpha\in\mathbb{R}$ are  closed.} and the
closedness properties are equivalent.\\

\noindent We say that $\varphi$ (resp. $\varphi^*$) is
differentiable  if it is differentiable  on
$]a_\varphi,b_\varphi[$ (resp. $]a_{\varphi^*},b_{\varphi^*}[$),
the interior of its domain. We say also that $\varphi$ (resp.
$\varphi^*$) is strictly convex   if it is strictly convex  on
$]a_\varphi,b_\varphi[$ (resp. $]a_{\varphi^*},b_{\varphi^*}[$).\\

\noindent The strict convexity of $\varphi$ is equivalent to the
condition that its conjugate $\varphi^*$ is essentially smooth,
i.e., differentiable with
\begin{equation}
\begin{array}{ccccc}
\lim_{t\downarrow a_{\varphi^*}}{\varphi^*}^\prime (t) & = &
-\infty &
\text{ if } & a_{\varphi^*}> -\infty,\\
\lim_{t\uparrow b_{\varphi^*}}{\varphi^*}^\prime (t) & = & +\infty
& \text{ if } & b_{\varphi^*}< +\infty.
\end{array}
\end{equation}
Conversely, $\varphi$ is essentially smooth  if and only if
$\varphi^*$ is strictly convex; see e.g. \cite{Rockafellar1970}
section 26 for the proofs of these properties.\\
\\
If $\varphi$ is differentiable, we denote $\varphi'$ the
derivative function of $\varphi$,  and we define
$\varphi'(a_\varphi)$ and $\varphi'(b_\varphi)$ to be the limits
(which may be finite or infinite) $\lim_{x\downarrow
a_\varphi}\varphi'(x)$ and $\lim_{x\uparrow
b_\varphi}\varphi'(x)$, respectively. We denote
$\text{Im}\varphi'$ the set of all values of the function
$\varphi'$, i.e., $\text{Im}\varphi':=\left\{\varphi'(x) \text{
such that } x\in [a_\varphi,b_\varphi] \right\}$. If additionally
the function $\varphi$ is strictly convex, then $\varphi'$ is
increasing on $[a_\varphi,b_\varphi]$. Hence, it is one-to-one
function from $[a_\varphi,b_\varphi]$ to $\text{Im}\varphi'$, we
denote in this case ${\varphi'}^{-1}$ the inverse function of
$\varphi'$ from
$\text{Im}\varphi'$ to $[a_\varphi,b_\varphi]$.\\

\noindent Note that if $\varphi$ is differentiable, then for all
$x\in ]a_\varphi,b_\varphi[$,
\begin{equation}\label{eqn 1}
\varphi^*\left(\varphi'(x)\right)=x\varphi'(x)-\varphi\left(
x\right).
\end{equation}
If additionally $\varphi$ is strictly convex, then for all $t\in
\text{Im}\varphi'$ we have
\begin{equation}
\varphi^*(t)=t{\varphi'}^{-1}(t)-\varphi\left({\varphi'}^{-1}(t)\right)\quad
\text{ and } \quad {\varphi^*}^\prime (t)={\varphi'}^{-1}(t).
\end{equation}
On the other hand, if $\varphi$ is essentially smooth, then the
interior of the domain of $\varphi^*$ coincides with that of
$\text{Im}\varphi'$, i.e.,
$\left(a_{\varphi^*},b_{\varphi^*}\right)=\left(\varphi'(a_\varphi),
\varphi'(b_\varphi)\right)$.\\

\noindent The domain of the $\phi$-divergence will be denoted
$\text{dom}\phi$, i.e.,
\begin{equation}
 \text{dom}\phi :=\left\{Q\in\mathcal{M} \text{ such that }\phi(Q,P)
 <\infty\right\}.
\end{equation}

\begin{definition}
Let $\Omega$ be some subset in $\mathcal{M}$. The
$\phi-$divergence between  the set $\Omega$ and  a p.m. $P$, noted
$\phi(\Omega,P)$,  is
\begin{equation*}
\phi(\Omega,P):=\inf_{Q\in\Omega}\phi(Q,P).
\end{equation*}
\end{definition}

\begin{definition}
Assume that $\phi(\Omega,P)$ is finite. A measure $Q^*\in\Omega$
such that
\begin{equation*}
    \phi(Q^*,P)\leq \phi(Q,P) ~\text{ for all }~Q\in\Omega
\end{equation*}
is called a $\phi$-projection of $P$ on $\Omega$. This projection
may not exist, or may be not defined uniquely.\\
\end{definition}

\noindent If $\varphi$ is a strictly convex, then the function
$Q\in\mathcal{M}(P)\mapsto\phi(Q,P)$ is strictly convex, and the
$\phi$-projection of $P$ on some convex set $\Omega$ is uniquely
defined whenever it exists.\\

\noindent Let $g_i : \mathcal{X} \mapsto \mathbb{R}$,
$i=1,\ldots,l$, be measurable real valued  functions on
$\mathcal{X}$. Denote  $g:= (g_0,g_1,\ldots,g_l)^T$ with
$g_0:=\mathds{1}_\mathcal{X}$. We assume that the functions $g_0,
g_1,\ldots,g_l$ are linearly independent in the following sense :
$P\left\{\lambda^Tg(x)\neq 0\right\}>0$ for any $\lambda\in
\mathbb{R}^{1+l}$ with $\lambda \neq 0$. For all
$\lambda\in\mathbb{R}^{1+l}$, we denote
$\lambda_0,\lambda_1,\ldots,\lambda_l$ the $(1+l)$ coordinates of
$\lambda$.\\

\noindent Let's denote by $M_g$ the set of all signed finite
measures with total mass one, a.c. w.r.t. $P$, which integrate the
functions $g_i$ and satisfy a finite number of linear constraints,
i.e.,
 \begin{equation}\label{M sets}
  M_g:=\left\{Q\in \mathcal{M}(P) \text{ such that }  Q(\mathcal{X})=1 \text{ and }
  \int_\mathcal{X}
  g_i(x)~dQ(x)=0, i=1,\ldots,l\right\}.
 \end{equation}
We consider the optimization problem
\begin{equation}\label{optim prob}
  \inf_{Q\in M_g}\phi(Q,P).
\end{equation}
The Lagrangian ``dual'' problem associated with (\ref{optim prob})
is
\begin{equation}\label{lagrangian dual prob}
\sup_{\lambda\in \mathbb{R}^{1+l}}\left\{ \lambda_0 -
\int_\mathcal{X}
\varphi^*\left(\lambda^Tg(x)\right)~dP(x)\right\}.
\end{equation}
\noindent We will consider the problem of the ``dual'' equality
$\inf (\ref{optim prob}) = \sup (\ref{lagrangian dual prob})$, the
existence of optimal solutions in (\ref{lagrangian dual prob}),
and in particular the problems of the existence and the
characterization of the optimal solutions in (\ref{optim prob}),
i.e.,
the $\phi$-projections of $P$ on the set $M_g$.\\
\\
These problems intervene in the domain of the reconstruction of
laws, in particular, the classical moment problem. Also they
appear frequently in Statistics; in fact, the recent ``empirical
likelihood'' method, which is the non parametric version of the
celebrate maximum likelihood method, can be seen as an
optimization problem of the $KL_m$-divergence between some set of
measures defined as in (\ref{M sets}) and
the empirical measure associated to the data.\\
\\
In the vocabulary of the duality theory, a measure $Q$ in $M_g$
which realizes the infimum in (\ref{optim prob}) (i.e., a
$\phi$-projection of $P$ on $M_g$ in the vocabulary of
$\phi$-divergences theory) is called ``a primal optimal solution''
or simply ``an optimal solution'', while a point $\lambda$ in
$\mathbb{R}^{1+l}$ realizing the supremum in (\ref{lagrangian
dual prob}) is called ``a dual optimal solution''.\\
\\
\noindent For the optimization problem of convex function $\psi :
\mathbb{R}^n \mapsto ]-\infty, +\infty]$ on convex sets $C$ in
$\mathbb{R}^n$ subject to linear constraints $Ax=b\in\mathbb{R}^m$
where $A$ is some $m\times n$-matrix, a sufficient condition, in
order that the equality
\begin{equation}\label{dual equality in Rn}
\inf_{\left\{x\in C;~Ax=b\right\}}\psi(x)=\sup_{t\in
\mathbb{R}^m}\left\{b^Tt-\psi^*\left(A^Tt\right)\right\}
\end{equation}
holds with dual attainment, is that there exists a point
$\widehat{x}$ in the relative interior\footnote{i.e., the interior
in the real affine subspace $\langle C\cap \text{dom}\psi\rangle$
of $\mathbb{R}^n$ endowed with the relative topology of the usual
topology on $\mathbb{R}^n$.}  of the convex set $C\cap
\text{dom}\psi$ such that $A\widehat{x}=b$. See e.g.
\cite{Rockafellar1970}  for the
proofs of these results.\\
\\
In order to make the set  $M_g$ closed and  the linear functions
$Q\in \mathcal{M} \mapsto \int_\mathcal{X}g_i(x)~dQ(x)$ continuous
(which we need to apply the duality theory and to treat the
problem of existence of $\phi$-projections of $P$ on the set
$M_g$), we endow the vector space $\mathcal{M}$ by the weak
topology which we denote $\tau_\mathcal{F}$ induced by
$\mathcal{F}\cup \mathcal{B}_b$ where
$\mathcal{F}:=\left\{g_0,g_1,\ldots,g_l\right\}$ and
$\mathcal{B}_b$ is the set of all bounded $\mathcal{B}$-measurably
real valued functions on $\mathcal{X}$; see section 2 below for
precise definition of the $\tau_{\mathcal{F}}$-topology.\\
\\
Note that the relative interior of the convex set $M_g$ is
generally empty in the weak topology $\tau_\mathcal{F}$.
\cite{BorweinLewis1992I} have extended the idea of the relative
interior (r.i.) of convex sets in $\mathbb{R}^n$ to a new notion
which have called ``the quasi relative interior'' (q.r.i.) of
convex subsets of an arbitrary Hausdorff topological vector  space
$X$ (having finite or infinite dimension), and  they used it to
construct a powerful duality theory for the optimization problem
of convex function $\psi : X \mapsto (-\infty,+\infty]$ on convex
sets $C\subseteq X$ subject to linear constraints. In particular,
when $X$ is locally convex, they obtain similar results as in
(\ref{dual equality in Rn}) when the relative interior is replaced
by the quasi relative interior; see \cite{BorweinLewis1992I}
Corollary 4.8. The main advantage of the quasi relative interior
of convex subset $C$ of infinite dimension vector space $X$ is
that it is frequently
nonempty even when the relative interior of $C$ is empty.\\
\\
If $\int_\mathcal{X} |g_i(x)|~dP(x)$ is finite for all
$i=1,\ldots,l$, then the convex conjugate of the convex function
$Q \mapsto \phi(Q,P)$ (on the vector space
$\mathcal{M}_\mathcal{F}(P)$ of all signed finite measures $Q$
a.c. w.r.t. $P$ and which integrate all the elements of
$\mathcal{F}$, i.e., all the functions $g_i$) can be written as
\begin{equation*}
\phi^*\left(f\right):=\sup_{Q\in\mathcal{M}_\mathcal{F}(P)}\left\{\int
f~dQ-\phi(Q,P)\right\}=\int \varphi^*\left(f\right)~dP, \text{ for
all } f\in \langle \mathcal{F}\cup \mathcal{B}_b\rangle;
\end{equation*}
see section 4 below for details. So, in this case, as in
\cite{BorweinLewis1991}, we can apply Corollary 4.8 of
\cite{BorweinLewis1992I} to obtain the dual equality $\inf
(\ref{optim prob}) = \sup (\ref{lagrangian dual prob})$ with dual
attainment, whenever there exists a measure $Q_0$ in $M_g$ which
belongs to the quasi relative interior of $M_g\cap
\text{dom}\phi$. This condition is called ``constraint
qualification''.  We can prove also from \cite{BorweinLewis1992I}
that a measure $Q_0$ is in the q.r.i of $M_g\cap \text{dom}\phi$
iff $a_\varphi < \frac{dQ_0}{dP}\leq \frac{dQ_0}{dP}<b_\varphi$,
$P$-almost everywhere ($P$-a.e.). We summarize these results and
some other results about the problems of the existence and the
characterization of the primal optimal solution (i.e., the
$\phi$-projection of $P$ on $M_g$) in the following two Theorems
and two Corollaries. For proofs, see Theorem 3.10 of
\cite{BorweinLewis1992I} , Corollary 2.6 and Theorem 4.8 of
\cite{BorweinLewis1991},  and Theorem II.2 of
\cite{CsiszarGamboaGassiat1999}.\\

\begin{theorem}\label{theorem BL}
If $\int_\mathcal{X} |g_i(x)|~dP(x)$ is finite for all
$i=1,\ldots,l$, and if the following constraint
qualification\footnote{The strict inequalities in (\ref{contraint
qualification}) mean that $P\left\{\frac{d\widehat{Q}}{dP}\leq
a_\varphi\right\}= P\left\{\frac{d\widehat{Q}}{dP}\geq
b_\varphi\right\}=0$.}:
 \begin{equation}\label{contraint qualification}
  \text{ there is a } \widehat{Q}\in M_g\cap\text{dom}\phi  \text{ such that }
  a_\varphi <  \frac{d\widehat{Q}}{dP}\leq
  \frac{d\widehat{Q}}{dP}< b_\varphi \quad (P-a.e.)
 \end{equation}
holds, then $\inf (\ref{optim prob}) = \sup (\ref{lagrangian dual
prob})$ and there is attainment in (\ref{lagrangian dual prob}).
Suppose additionally that $\varphi^*$ is essentially smooth (which
is equivalent to the strict convexity of $\varphi$), and that
there exists a dual optimal solution $\overline{\lambda}$ which is
an interior point of
\begin{equation}\label{dom de phi star dans R l +1}
\text{dom} \phi^* :=\left\{\lambda\in\mathbb{R}^{1+l}~\text{ such
that }\int_\mathcal{X} \varphi^*\left(\lambda^Tg(x)\right)~dP(x)
\text{ is finite }\right\}.
\end{equation}
Then the unique optimal solution of (\ref{optim prob}) (i.e., the
$\phi$-projection of $P$ on $M_g$), which we denote by $Q^*$,
exists and it is given by
\begin{equation}\label{characterization}
 \frac{dQ^*}{dP}(x) =
 {\varphi ^*}^\prime\left(\overline{\lambda}^Tg(x)\right).\\
\end{equation}
\end{theorem}

\noindent In (\ref{dom de phi star dans R l +1}), for brevity, the
definition of $\text{dom}\phi^*$, which usually is the set of
functions $f$ such that $\phi^*(f)<\infty$, is modified here.\\

\begin{remark}
If all functions  $g_i$ belong to $L_\infty(\mathcal{X},P)$, and
if for a dual optimal solution
$\overline{\lambda}\in\mathbb{R}^{1+l}$ the following condition
\begin{equation}\label{integrability sufficient condition 1}
a_{\varphi^*}< \text{ess}\inf \overline{\lambda}^Tg(.)\leq
\text{ess}\sup \overline{\lambda}^Tg(.) < b_{\varphi^*}
\end{equation}
holds, then $\overline{\lambda}$ is an interior point of
$\text{dom}\phi^*$. Hence, under assumption (\ref{integrability
sufficient condition 1}), all results in the above Theorem hold
whenever the constraint qualification (\ref{contraint
qualification}) is met.\\
\end{remark}

\noindent If all functions  $g_i$ belong to
$L_\infty(\mathcal{X},P)$, and the convex function $\varphi^*$ is
everywhere finite (i.e., $a_{\varphi^*}=-\infty$ and
$b_{\varphi^*}=+\infty$), then obviously condition
(\ref{integrability sufficient condition 1}) holds since
$\text{dom}\phi^*=\mathbb{R}^{1+l}$ in this case. Hence, under the
constraint qualification (\ref{contraint qualification}), all
results in the above Theorem hold. We state this result in the
following Corollary.\\

\begin{corollary}\label{corollary theorem BL}
Suppose that all functions  $g_i$ belong to
$L_\infty(\mathcal{X},P)$ and  $\varphi^*$ is everywhere finite
(i.e., $a_{\varphi^*}=-\infty$ and $b_{\varphi^*}=+\infty$). If
the constraint qualification (\ref{contraint qualification})
holds, then $\inf (\ref{optim prob}) = \sup (\ref{lagrangian dual
prob})$ and there is attainment in (\ref{lagrangian dual prob}).
Suppose additionally that $\varphi^*$ is everywhere differentiable
(which is equivalent to the strict convexity of $\varphi$), then
the unique optimal solution $Q^*$ of (\ref{optim prob}) (i.e., the
$\phi$-projection of $P$ on $M_g$)  exists and it is given by
\begin{equation}\label{characterization}
 \frac{dQ^*}{dP}(x) =
 {\varphi^*}^\prime\left(\overline{\lambda}^Tg(x)\right),
\end{equation}
where $\overline{\lambda}\in\mathbb{R}^{1+l}$ is any dual optimal
solution.
\end{corollary}

\noindent In the following Theorem and Corollary, we give
sufficient conditions for the uniqueness of the dual optimal
solution (see \cite{BorweinLewis1991} Theorem 4.5 for the proof).
Note that the strict convexity of $\varphi^*$ is equivalent to the
condition that its conjugate $\varphi$ is essentially smooth.\\

\begin{theorem}\label{theorem BL uniqueness 1}
Suppose that all assumptions of Theorem \ref{theorem BL} are
satisfied. Suppose furthermore that the function $\varphi$ is
essentially smooth. Then the dual optimal solution
$\overline{\lambda}$ is unique. Moreover, the unique optimal
solution  $Q^*$ of (\ref{optim prob})
 exists and it is given by
\begin{equation}\label{characterization 2}
 \frac{dQ^*}{dP}(x) =
 {\varphi ^*}^\prime\left(\overline{\lambda}^Tg(x)\right)={\varphi'}^{-1}\left(
 \overline{\lambda}^Tg(x)\right).
\end{equation}
\end{theorem}

\begin{corollary}\label{theorem BL uniqueness 1}
Suppose that all assumptions of Corollary \ref{corollary theorem
BL} are satisfied. Suppose additionally that the function
$\varphi$ is essentially smooth. Then the dual optimal solution
$\overline{\lambda}$ is unique. Moreover, the unique optimal
solution  $Q^*$ of (\ref{optim prob}) exists and it is given by
\begin{equation}\label{characterization 2}
\frac{dQ^*}{dP}(x) =
 {\varphi ^*}^\prime\left(\overline{\lambda}^Tg(x)\right)={\varphi'}^{-1}\left(
 \overline{\lambda}^Tg(x)\right).
\end{equation}\\
\end{corollary}
\noindent The important Corollary \ref{corollary theorem BL},
which essentially requires  that the constraint qualification
(\ref{contraint qualification}) holds, applies in the
$KL$-divergence case since the corresponding conjugate $\varphi^*$
is everywhere finite  (see also \cite{BorweinLewis1993} for other
examples), but it fails in the two important cases of Burg
relative entropy ($KL_m$-divergence in the context of divergences)
and Hellinger divergence without additional conditions since the
corresponding conjugates $\varphi^*$ are infinite on the intervals
$[1,+\infty)$ and $[2,+\infty)$, respectively.\\
\\
\noindent \cite{Leonard2001c} consider the optimization problem
(\ref{optim prob}) when the set $M_g$ is replaced by the subset
\begin{equation}\label{Mos}
M_{os}:=\left\{Q\in\mathcal{M}(P)\text{ such that }
q:=\frac{dQ}{dP}\in L_{\varphi^{**}_m}, ~
\int_\mathcal{X}g(x)~dQ=(1,0,\ldots,0)^T\right\},
\end{equation}
where $L_{\varphi^{**}_m}$ is the Orlicz space defined as follows:
\begin{equation}\label{the Orlicz space L}
L_{\varphi^{**}_m}:=\left\{q:\mathcal{X}\to \mathbb{R}; \text{
measurable such that } \|q\|_{\varphi^{**}_m}<\infty\right\}
\end{equation}
\begin{equation*}
\text{ with } ~~ \|q\|_{\varphi^{**}_m}:=\inf\left\{a>0;
\int_\mathcal{X}\varphi^{**}_m\left(\frac{|q(x)|}{a}\right)~dP(x)\leq
1 \right\},
\end{equation*}
and $\varphi^{**}_m$ is the convex conjugate of the convex
function $\varphi^*_m$ defined by
$\varphi^*_m(t):=\max\left(\varphi^*(t),\varphi^*(-t)\right)$ for
all $t\in\mathbb{R}$. Without the constraint qualification
(\ref{contraint qualification}), under the following integrability
condition
\begin{equation}\label{integrability condition}
\text{ for any } \lambda\in\mathbb{R}^{1+l}, ~
\int_\mathcal{X}\varphi^*\left(\lambda^Tg(x)\right)~dP(x)<\infty,
\end{equation}
applying the duality theory on Orlicz spaces, \cite{Leonard2001c}
obtains the dual equality
\begin{equation}\label{dual equality os}
\inf_{Q\in
M_{os}}\phi(Q,P)=\sup_{\lambda\in\mathbb{R}^{1+l}}\left\{\lambda_0-
\int_\mathcal{X}\varphi^*\left(\lambda^Tg(x)\right)~dP(x)\right\}.
\end{equation}
Moreover, if the value is finite, then  there exists at least one
$\phi$-projection of $P$ on $M_{os}$; see Theorem 3.4 of
\cite{Leonard2001c} for details in more general context. A
characterization of the $\phi$-projections of $P$ on the set
$M_{os}$ (with finite or infinite number of linear constraints) is
stated by \cite{Leonard2001b} under condition (\ref{integrability
condition}); see Theorems 4.4, 4.5 and 4.6 of \cite{Leonard2001b}.
Note that the integrability condition (\ref{integrability
condition}) implies that $\varphi^*$ is everywhere finite, and
these results apply in the important $KL$-divergence case with
finite or infinite number of linear constraints. However, the
condition (\ref{integrability condition}) does not hold in the
$KL_m$-divergence and Hellinger divergence cases since the domains
of the corresponding $\varphi^*$ functions are proper subsets of
$\mathbb{R}$, and the important result (\ref{dual equality os})
does not apply in these two important cases. Under the weaker
integrability assumption
\begin{eqnarray}
\text{ for any } \lambda\in\mathbb{R}^{1+l}, ~\text{ there exists
} \alpha >0 ~\text{ such that }~\nonumber\\
 \int_\mathcal{X}\varphi^*\left(\alpha\lambda^Tg(x)\right)~dP(x)
 +\int_\mathcal{X}\varphi^*\left(-\alpha\lambda^Tg(x)\right)~dP(x)<\infty,
 \label{weaker integrability condition}
\end{eqnarray}
the dual equality (\ref{dual equality os}) may fail;
see Theorem 3.3 of \cite{Leonard2001c}.\\
\\
\noindent The goal of the present paper is to give  results of
existence and characterization of the $\phi$-projections of a
given p.m. $P$ on some subsets $\Omega$ of
$\mathcal{M}_\mathcal{F}$, the space of all signed finite measures
which integrate a given class $\mathcal{F}$ of functions, in
particular, convex sets of signed finite measures defined by
linear constraints as in (\ref{M sets}) extending some previous
works (about the existence and characterization of the
$\phi$-projections on subsets of $\mathcal{M}^1$, the set of all
p.m.'s) of \cite{Csiszar1975}, \cite{Liese1977},
\cite{Csiszar1984}, \cite{Ruschendorf1984},
\cite{Ruschendorf1987}, \cite{Liese-Vajda1987},
\cite{TeboulleVajda1993} and \cite{Csiszar1995}. We give also
different versions of dual representations of the
$\phi$-divergences viewed as convex functions on the space of all
signed finite measures which integrate an arbitrary class of
functions. When the set $\Omega$ is defined by linear constraints
as in (\ref{M sets}), we consider the dual problem, and we obtain
the equality $\inf (\ref{optim prob}) = \sup (\ref{lagrangian dual
prob})$ with dual attainment under different assumptions without
constraint qualification. Additional conditions are given to
obtain similar results which apply in the two important
$KL_m$-divergence and Hellinger
divergence cases.\\

\noindent Enhancing $\mathcal{M}^1$ to $\mathcal{M}$ is motivated
by the following arguments: sometimes the $\phi$-projection, say
$Q^*_1$, of a p.m. $P$ on subset of $\mathcal{M}^1$ is not an
``interior'' point and we can not give in this case a definite
description of $Q^*_1$, while the $\phi$-projection, say $Q^*$, of
a p.m. $P$ on subset of $\mathcal{M}$ is an ``interior'' point,
which allows to give a perfect characterization of the
$\phi$-projection $Q^*$  (see  example \ref{contre exemple}).  In
the context of statistical estimation and tests using the
empirical likelihood method (see \cite{Owen2001}), or related ones
to criterions defined through divergences (see
\cite{Broniatowski-Keziou2004}), the projection of the empirical
measure $P_n$ of a sample on a set $\Omega^1$ of p.m.'s may make
problems when the projection is not an interior point of
$\Omega^1\cap\text{dom}\phi(.,P_n)$. Enhancing $\mathcal{M}^1$ to
$\mathcal{M}$, this difficulty does not hold any longer, and tests
as well as estimation can be performed.\\
\\
\noindent The rest of this paper is organized as follows : In
section 2, we consider the problem of existence of
$\phi$-projections on general closed sets of signed measures. In
section 3, we deal with the problem of characterization of the
$\phi$-projections on sets of signed measures, in particular, sets
of signed measures defined by linear constraints. In section 4, we
give different dual representations of $\phi$-divergences seen as
convex functions on the vector space of all signed finite measures
which integrate a given class of functions. In section 5, we apply
the results of sections 2, 3 and 4, to obtain the dual equality
$\inf (\ref{optim prob}) = \sup (\ref{lagrangian dual prob})$ with
dual attainment, under different assumptions without constraint
qualification.
\section{Existence of  $\phi$-Projections on Sets
of signed measures} \noindent In this section, we give sufficient
conditions for the existence of $\phi$-projections of some p.m.
$P$ on  sets $\Omega$ of signed finite measures which integrate
some class of functions (see Theorems \ref{corollary thm exist
proj sur ferme},
 \ref{thm exist proj sur ferme 1} and \ref{thm exist proj sur ferme},
 and Corollary \ref{existence
de proj pour klm and h} below). At first, we give some notation
and we establish a convenient topological context for this
problem. Let $\mathcal{F}$ be some class of measurable real valued
functions $f$ (bounded or unbounded) defined on $\mathcal{X}$.
Here, $\mathcal{F}$ is not assumed to be finite. Denote by
$\mathcal{B}_b$ the set of all bounded measurable real valued
functions defined on $\mathcal{X}$, and by
$\langle\mathcal{F}\cup\mathcal{B}_b\rangle$ the linear span of
$\mathcal{F}\cup \mathcal{B}_b$. Define the set
\begin{equation*}
\mathcal{M}_{\mathcal{F}}^{1}:=\left\{ Q\in \mathcal{M}^{1}~\text{
such that }~\int |f|~dQ<\infty ,\text{ for all  }f\text{  in
}\mathcal{F} \right\},
\end{equation*}
and the real vector space
\begin{equation*}
\mathcal{M}_{\mathcal{F}}:=\left\{ Q\in \mathcal{M}~\text{ such
that  }~\int |f|~d|Q|<\infty ,\text{  for all  }f\text{  in
}\mathcal{F}\right\},
\end{equation*}
in which $|Q|$ denotes the total variation of the signed finite
measure $Q$.\\
Note that if  $\mathcal{F}=\mathcal{B}_{b}$, then
$\mathcal{M}_{\mathcal{F}}^{1}=\mathcal{M}^{1}$ and
$\mathcal{M}_{\mathcal{F}}=\mathcal{M}$.
\begin{definition}\label{def tau phi}
Denote by $\tau_{\mathcal{F}}$ the weakest topology on
$\mathcal{M}_{\mathcal{F}}$ for which all mappings $Q\in
\mathcal{M}_{\mathcal{F}}\mapsto \int f~dQ$ are continuous when
$f$ belongs to $\mathcal{F}\cup\mathcal{B}_b$. Denote also by
$\tau_\mathcal{M}$ the weakest topology on
$\langle\mathcal{F}\cup\mathcal{B}_b\rangle$ for which all
mappings $f\in \langle\mathcal{F}\cup\mathcal{B}_b\rangle \mapsto
\int f~dQ$ are continuous when $Q\in\mathcal{M}_\mathcal{F}$. We
sometimes call $\tau_\mathcal{F}$ the topology induced  by
$\langle\mathcal{F}\cup\mathcal{B}_b\rangle$ on
$\mathcal{M}_\mathcal{F}$, and likewise $\tau_\mathcal{M}$ the
topology induced by $\mathcal{M}_\mathcal{F}$ on
$\langle\mathcal{F}\cup\mathcal{B}_b\rangle$.
\end{definition}

\noindent A base of open neighborhoods for any $R$ in
$\mathcal{M}_{\mathcal{F}}$ is defined by
\begin{equation}\label{base de voisinages}
U(R,\mathcal{A},\varepsilon):=\left\{Q\in
\mathcal{M}_{\mathcal{F}}~\text{ such that }~\max_{f\in
\mathcal{A}}\left|\int f~dR-\int f~dQ\right| <\varepsilon \right\}
\end{equation}
for $\varepsilon >0$ and $\mathcal{A}$ a finite collection of
functions in $\langle\mathcal{F}\cup\mathcal{B}_b\rangle$.\\
\\
\noindent We refer to  Chapter 5 of \cite{Dunford-Schwartz1962},
for the various topologies induced by classes of functions. Note
that the class $\mathcal{B}_b$ induces the so-called
$\tau$-topology (see e.g.
\cite{Groeneboom-Oosterhoff-Ruymgaart1979} and
\cite{Ganssler1971}), and that $\mathcal{M}_{\mathcal{B}_b}$ is
the whole vector space $\mathcal{M}$.\\
\\
\noindent The above $\tau_{\mathcal{F}}-$topology on
$\mathcal{M}_{\mathcal{F}}$ is indeed the natural  and the most
convenient one in order to handle projection properties. It has
been introduced in the context of large deviation probabilities by
\cite{EischelsbacherScmock1997} for the Kullback-Leibler
divergence and it is used  in Statistics in
\cite{Broniatowski2002}, \cite{Keziou2003},
\cite{Broniatowski-Keziou2003} and \cite{Keziou2003b}. Usually the
sets which are to be considered in statistical applications are
not compact but merely closed sets; a typical example is when they
are defined by linear constraints as in (\ref{M sets}). Hence, the
set $M_g$ is closed in $\mathcal{M}_\mathcal{F}$ endowed with the
$\tau_\mathcal{F}$-topology if the  functions $g_i$ (which may be
bounded or unbounded) belong to $\mathcal{F}$; this motivates the
choice of $\tau_\mathcal{F}$-topology.\\

\begin{proposition}\label{evthlc}
Equip $\mathcal{M}_{\mathcal{F}}$ with the
$\tau_{\mathcal{F}}$-topology and
$\langle\mathcal{F}\cup\mathcal{B}_b\rangle$ with the
$\tau_\mathcal{M}$-topology. Then, $\mathcal{M}_{\mathcal{F}}$ and
$\langle\mathcal{F}\cup\mathcal{B}_b\rangle$  are Hausdorff
locally convex topological vector spaces. Further, the topological
dual space of $\mathcal{M}_{\mathcal{F}}$  is the set of all
mappings $Q\mapsto \int f~dQ$ when $f$ belongs to
$\langle\mathcal{F}\cup\mathcal{B}_b\rangle$, and the topological
dual space of $\langle\mathcal{F}\cup\mathcal{B}_b\rangle$   is
the set of all mappings $f\mapsto \int f~dQ$ when $Q$ belongs to
$\mathcal{M}_{\mathcal{F}}$.\\
\end{proposition}

\noindent  \textbf{Proof of Proposition \ref{evthlc}} By Lemma
5.3.3 in \cite{Dunford-Schwartz1962}, the vector space
$\mathcal{M}_{\mathcal{F}}$ equipped with the
$\tau_{\mathcal{F}}$-topology is a Hausdorff locally convex
topological space. On the other hand, the set of all mappings
$Q\in\mathcal{M}_\mathcal{F}\mapsto \int f~dQ$ when $f$  belongs
to $\langle\mathcal{F}\cup \mathcal{B}_b\rangle$ is a total linear
space; indeed, for all $Q\in\mathcal{M}_\mathcal{F}$, assume that
$\int f~dQ=0$ for all $f$ in $<\mathcal{F}\cup\mathcal{B}_b>$,
choose $f=\mathds{1}_{\left\{B\right\}}$ for any $B\in\mathcal{B}$
to conclude that $Q=0$. The proof ends then as a consequence of
Theorem 5.3.9 in \cite{Dunford-Schwartz1962}. $\blacksquare$\\
\\
\noindent  We denote by
$\left[\mathcal{M}_\mathcal{F};\tau_\mathcal{F}\right]$ and by
$\left[\langle\mathcal{F}\cup\mathcal{B}_b\rangle;
\tau_\mathcal{M}\right]$ the two Hausdorff locally convex
topological vector spaces endowed with the
$\tau_\mathcal{F}$-topology and the $\tau_\mathcal{M}$-topology,
respectively.\\
\\
\noindent \cite{Broniatowski-Keziou2003} have proved that the
function
$Q\in\left[\mathcal{M}_\mathcal{F};\tau_\mathcal{F}\right]\mapsto
\phi(Q,P)$ is lower semi-continuous (l.s.c.), provided only that
the corresponding convex function $\varphi$ is closed; see
Proposition 2.3 of  \cite{Broniatowski-Keziou2003}  and
Proposition 2.1 of  \cite{Keziou2003}  which we recall here for
convenience.\\

\begin{proposition}\label{prop div lsc param}
For any $\phi$-divergence, the divergence function $Q\mapsto
\phi(Q,P)$ from $\left[\mathcal{M}_{\mathcal{F}};\tau
_{\mathcal{F}}\right]$ to $[0,+\infty]$ is l.s.c.\\
\end{proposition}

\noindent We will use the following Lemma to prove Proposition
\ref{prop div lsc param}.\\

\begin{lemma}\label{Lemma M-P-F ferme}
Let $\mathcal{M}_{\mathcal{F}}(P)$ denotes the vector subspace of
all signed measures in $\mathcal{M}_{\mathcal{F}}$ which are
absolutely continuous w.r.t. $P$. The vector subspace
$\mathcal{M}_{\mathcal{F}}(P)$ is a closed set in
$\left[\mathcal{M}_{\mathcal{F}};\tau_{\mathcal{F}}\right]$.\\
\end{lemma}
\noindent \textbf{Proof of Lemma \ref{Lemma M-P-F ferme}} Let
$\overline{\mathcal{M}_{\mathcal{F}}(P)}$ denotes the closure of
$\mathcal{M}_{\mathcal{F}}(P)$ in
$\left[\mathcal{M}_{\mathcal{F}};\tau_{\mathcal{F}}\right]$.
Assume that there exists $R$ in
$\overline{\mathcal{M}_{\mathcal{F}}(P)}$ with $R$ not in
$\mathcal{M}_{\mathcal{F}}(P)$. Then, there exists some $B$ in
$\mathcal{B}$ such that  $P(B)=0$ and $R(B)\neq 0$. On the other
hand, for all $n$ in $\mathbb{N}$, the set
$U:=U\left(R,\mathds{1}_{\left\{B\right\}},1/n\right)$ is a
neighborhood of $R$ (see (\ref{base de voisinages})),  hence,
$U\cap \mathcal{M}_{\mathcal{F}}(P)$ is non void. Therefore, we
can construct a sequence of measures $R_{n}$ in
$\mathcal{M}_{\mathcal{F}}(P)$ such that
\begin{equation*}
\left|\int \mathds{1}_{\left\{B\right\}}~dR-\int
\mathds{1}_{\left\{B\right\}}~dR_{n}\right|<1/n.
\end{equation*}
Since $R_{n}(B)=0$ for all $n$ in $\mathbb{N}$, we deduce that
$R(B)=0$, a contradiction. This implies that
$\overline{\mathcal{M}_{\mathcal{F}}(P)}=\mathcal{M}_{\mathcal{F}}(P)$,
that is $\mathcal{M}_{\mathcal{F}}(P)$ is closed in
$\left[\mathcal{M}_{\mathcal{F}};\tau_\mathcal{F}\right]$. This
concludes the proof of Lemma \ref{Lemma M-P-F ferme}.
$\blacksquare$\\

\begin{remark}\label{remarque fermeture}
Note that if $\mathcal{F}=\mathcal{B}_b$, then
$\mathcal{M}_\mathcal{F}=\mathcal{M}$ and
$\mathcal{M}_\mathcal{F}(P)=\mathcal{M}(P)$. Hence, we deduce from
Lemma \ref{Lemma M-P-F ferme} that the subspace $\mathcal{M}(P)$
is closed in $\left[\mathcal{M};\tau\right]$, the space of all
signed finite measures endowed with the $\tau$-topology. Note also
that $\mathcal{M}^1_\mathcal{F}$ and
$\mathcal{M}^1_{\mathcal{F}}(P)$ are closed in
$\left[\mathcal{M}_\mathcal{F};\tau_{\mathcal{F}}\right]$, and
that $\mathcal{M}^1$ and $\mathcal{M}^1(P)$ are closed in
$\left[\mathcal{M};\tau\right]$.\\
\end{remark}

\noindent \textbf{Proof of Proposition \ref{prop div lsc param}}
Let $\alpha$ be a real number. We  prove that the set
\begin{equation*}
A(\alpha):=\left\{Q\in \mathcal{M}_{\mathcal{F}}~\text{ such that
}~\phi(Q,P)\leq\alpha\right\}
\end{equation*}
is closed in $\left[\mathcal{M}_{\mathcal{F}};\tau
_{\mathcal{F}}\right]$. By Lemma \ref{Lemma M-P-F ferme},
$\mathcal{M}_{\mathcal{F}}(P)$ is closed in $\left[
\mathcal{M}_{\mathcal{F}};\tau_{\mathcal{F}}\right]$. Since
$A(\alpha)$ is included in
$\left[\mathcal{M}_{\mathcal{F}}(P);\tau_{\mathcal{F}}\right]$, we
have to prove that $A(\alpha)$ is closed in the subspace
$\left[\mathcal{M}_{\mathcal{F}}(P);\tau_{\mathcal{F}}\right]$.
Let
\begin{equation*}
B(\alpha):=\left\{f\in L_{1}(\mathcal{X},P)~\text{  such that
}~\int \varphi(f(x))~dP(x)\leq\alpha\right\}.
\end{equation*}
$B(\alpha)$ is a convex set, since $\varphi$ is a convex function.
Furthermore, $B(\alpha)$ is closed in $L_{1}(\mathcal{X},P)$.
Indeed, let $f_{n}$ be a sequence in $B(\alpha)$ with
$\lim_{n\rightarrow\infty}f_{n}=f^*$, where the limit is intended
in $L_{1}(\mathcal{X},P)$. Hence, there exists a subsequence
$f_{n_{k}}$ which converges to $f^*$ ($P$-a.e.). The functions
$\varphi(f_{n_{k}})$ are nonnegative. Further, we have
$\liminf_{k\rightarrow +\infty }\varphi (f_{n_{k}}(x))=f^*(x)$~
($P$-a.e.) by the closedness of the convex function $\varphi$.
Therefore, Fatou's Lemma implies
\begin{equation*}
\int \varphi(f^*)~dP~\leq \int \liminf_{k\rightarrow +\infty
}\varphi (f_{n_{k}})~dP~\leq \liminf_{k\rightarrow +\infty }\int
\varphi (f_{n_{k}})~dP~\leq \alpha,
\end{equation*}
which is to say that $f^*$  belongs to $B(\alpha)$. Hence,
$B(\alpha)$ is a closed subset in $L_{1}(\mathcal{X},P)$. Since,
it is convex, it is then weakly closed in $L_{1}(\mathcal{X},P)$;
see e.g.  Theorem 5.3.13 in \cite{Dunford-Schwartz1962}. Denote by
$W$ the weak topology on $L_{1}(\mathcal{X},P)$ and consider the
mapping $H$ defined by
\begin{equation*}
\begin{array}{ccccc}
H & : &
\left[\mathcal{M}_{\mathcal{F}}(P);\tau_{\mathcal{F}}\right] &
\mapsto
& \left[L_{1}(\mathcal{X},P);W\right] \\
~ & ~ & Q & \mapsto & H(Q)=dQ/dP.
\end{array}
\end{equation*}
Let us prove that $H$ is weakly continuous, that is $~Q\mapsto
\int H(Q)g~dP~$ is a  continuous mapping for all $g$ in
$L_\infty(\mathcal{X},P)$. Indeed, let $g$  be some function in
$L_\infty(\mathcal{X},P)$. Then, we have
\begin{equation*}
\int H(Q)g~dP=\int (dQ/dP)g~dP=\int g~dQ.
\end{equation*}
The mapping $Q\mapsto \int g~dQ$ is
$\tau_{\mathcal{F}}$-continuous; indeed, for all $g$ in
$L_\infty(\mathcal{X},P)$, it holds $P(g>{\Vert
g\Vert}_{\infty})=0$, which implies $Q(g>{\Vert
g\Vert}_{\infty})=0$, for all $Q$ in
$\mathcal{M}_{\mathcal{F}}(P)$. Therefore, $\int g~dQ=\int
g\mathds{1}_{[g\leq {\Vert g\Vert}_{\infty}]}~dQ$. Now, the
mapping $Q\mapsto \int g\mathds{1}_{[g\leq {\Vert
g\Vert}_{\infty}]}~dQ$ is continuous in
$\tau_\mathcal{F}$-topology since $g1_{[g\leq {\Vert
g\Vert}_{\infty}]}$$\in$$\mathcal{F}\cup\mathcal{B}_b$.  Since
$A(\alpha)=\left\{ Q\in
\mathcal{M}_{\mathcal{F}}(P),~\phi(Q,P)\leq\alpha\right\}=H^{-1}\left(
B(\alpha)\right)$, we deduce that $A(\alpha)$  is closed in
$\left[\mathcal{M}_{\mathcal{F}}(P);\tau_{\mathcal{F}}\right]$,
for any $\alpha$ in $\mathbb{R}$. This proves Proposition
\ref{prop div lsc
param}.$\blacksquare$\\
\\
\noindent For any $\phi$-divergence, by the lower semi-continuity
of the function
$Q\in\left[\mathcal{M}_\mathcal{F};\tau_\mathcal{F}\right]\mapsto
\phi(Q,P)$, the following result holds.\\

\begin{theorem}\label{prop proj sur compact}
Let $P$ be some p.m. and $\Omega$ some compact subset of
$\left[\mathcal{M}_{\mathcal{F}};\tau _{\mathcal{F}}\right]$. Then
there exists at least one $\phi$-projection of $P$ on $\Omega$.\\
\end{theorem}

\noindent Using some similar arguments as used in the proof of
Theorem 2.4 in \cite{Liese1977} or Proposition 8.5 in
\cite{Liese-Vajda1987} and Fenchel's inequality or Hölder
inequality, we state  general results for the existence of
$\phi$-projections of some p.m. $P$ on closed sets $\Omega$ of
$\left[\mathcal{M}_\mathcal{F};\tau_\mathcal{F}\right]$ (see
Theorem \ref{thm exist proj sur ferme 1} and \ref{thm exist proj
sur ferme} below). At first, in the following Theorem, we give a
version of Theorem 2.4 in \cite{Liese1977} or Proposition 8.5 in
\cite{Liese-Vajda1987}.\\

\begin{theorem}\label{corollary thm exist proj
sur ferme} Let $\Omega$ be some closed set in
$\left[\mathcal{M};\tau\right]$. Assume that the following two
conditions
\begin{equation}\label{phi Omega P est fini corollaire}
\phi(\Omega,P):=\inf_{Q\in\Omega}\phi(Q,P) <\infty\footnote{Note
that this is equivalent to the condition:
 there exists $Q\in \Omega$  such that
$\phi(Q,P)$ is finite.}
\end{equation}
and
\begin{equation}\label{condition pr uniforme integr corollaire}
\lim_{\left\vert x\right\vert \rightarrow \infty
}\frac{\varphi(x)}{\left\vert x\right\vert }= +\infty
\end{equation}
hold. Then there exists at least one $\phi$-projection of $P$ on
$\Omega$.\\
\end{theorem}

\noindent \textbf{Proof of Theorem \ref{corollary thm exist proj
sur ferme}} Denote $m:=\phi(\Omega,P)$ which is finite by
assumption, and let $\beta$ be a positive number. Define the sets
\begin{equation*}
\Omega (\beta ):=\left\{ Q\in\Omega ~\text{ such that
}~\phi(Q,P)\leq m+\beta \right\}
\end{equation*}
and
\begin{equation*}
\Lambda(\beta ):=\left\{ q:=\frac{dQ}{dP}~\text{ such that }~Q\in
\Omega(\beta)\right\}.
\end{equation*}
The set $\Lambda(\beta)$ is uniformly integrable. Hence, it is
weakly sequentially compact in $L_1(\mathcal{X},P)$, (see e.g.
\cite{Meyer1966} p. 39). Consider now a sequence $Q_{n}$ in
$\Omega(\beta)$ such that
\begin{equation*}
\lim_{n\rightarrow +\infty}\phi(Q_{n},P)=\phi(\Omega,P).
\end{equation*}
The sequence $q_{n}:=dQ_{n}/dP$ belongs to $\Lambda(\beta)$.
Therefore, there exists a subsequence
$\left(q_{n_i}\right)_{i\in\mathbb{N}}$ which converges weakly in
$L_1(\mathcal{X},P)$ to some function $q^{\ast}\in
L_1(\mathcal{X},P)$, which is to say that the corresponding
sequence of signed finite measures $Q_{n_i}$ converges to
$Q^{\ast}\in\mathcal{M}(P)$ in $\tau$-topology where $Q^{\ast}$ is
defined by $dQ^{\ast}/dP:=q^{\ast}$. Hence, $Q^*$ belongs to
$\Omega$ since it is the limit in $\tau$-topology of the sequence
$(Q_{n_i})$ which belongs to the $\tau$-closed set $\Omega$. On
the other hand, the mapping $Q\in\left[\mathcal{M};\tau\right]
\mapsto\phi(Q,P)$ is l.s.c.\footnote{this holds from  Proposition
\ref{prop div lsc param} choosing the class of functions
$\mathcal{F}=\mathcal{B}_b$, the class of all bounded measurable
real valued functions.}, and therefore
\begin{equation}\label{inequality by lsc}
\phi(Q^{\ast},P)\leq\lim_{i\to
+\infty}\phi(Q_{n_i},P)=\phi(\Omega,P)<\infty.
\end{equation}
We deduce that  $Q^*$ is a $\phi$-projection of $P$ on $\Omega$.
$\blacksquare$\\

\begin{remark}
For sets $\Omega$ of p.m.'s defined by linear constraints,
sufficient conditions for the existence of $KL$-projections are
presented in (\cite{Csiszar1975} Theorem 3.1,  Corollary 3.1  and
Theorem 3.3). Sufficient conditions of the existence of
$\phi$-projections on sets of p.m.'s satisfying linear equality or
inequality constraints are given in \cite{Csiszar1995} Theorem
3.\\
\end{remark}

\begin{remark}\label{ES result}
By \cite{EischelsbacherScmock1997}, if for all $\alpha
>0$ and all $f \in \mathcal{F}$,  $\int \exp
\left({\alpha|f|}\right)~dP<\infty$, then the level sets
\begin{equation*}
\left\{ Q\in \mathcal{M}_{\mathcal{F}}^{1}~\text{ such that }
~KL(Q,P)\leq c\right\}
\end{equation*} are compact in
$\left[\mathcal{M}_{\mathcal{F}}^{1};\tau_{\mathcal{F}}\right]$
for all real $c$. Therefore, for any $\tau_\mathcal{F}$-closed set
$\Omega\subset \mathcal{M}_{\mathcal{F}}^{1}$  for which
$KL(\Omega ,P)<\infty$, the projection of $P$ on $\Omega$ exists;
see \cite{EischelsbacherScmock1997} Lemma 2.1.\\
\end{remark}

\noindent Using Fenchel's inequality and some similar arguments to
that in Lemma 2.1 of \cite{EischelsbacherScmock1997},  We
generalize  Theorem 3 of \cite{Csiszar1995} and  the result in
Remark \ref{ES result} about the existence of projections, to the
class of $\phi$-divergences and to $\tau_{\mathcal{F}}$-closed
sets of signed measures.\\

\begin{theorem}\label{thm exist proj sur ferme 1}
Let $\Omega$ be some closed set in $\mathcal{M}_{\mathcal{F}}$
equipped with the $\tau _{\mathcal{F}}$-topology. Suppose that the
following three assumptions
\begin{equation}\label{phi Omega P est fini 1}
\phi(\Omega,P)<\infty,
\end{equation}
\begin{equation}\label{condition pr uniforme integr}
\lim_{\left\vert x\right\vert \rightarrow \infty
}\frac{\varphi(x)}{\left\vert x\right\vert }= +\infty
\end{equation}
\begin{equation}\label{condition Fenchel ineq}
\text{ and for every } f\in\mathcal{F} \text{ and every }
\alpha>0, \int \varphi^*\left(\alpha|f|\right)~dP<\infty
\end{equation}
hold. Then there exists at least one $\phi$-projection of $P$ on
$\Omega$.\\
\end{theorem}

\noindent \textbf{Proof of Theorem \ref{thm exist proj sur ferme
1}} As shown in the proof of Theorem \ref{corollary thm exist proj
sur ferme}, under assumptions (\ref{phi Omega P est fini 1}) and
(\ref{condition pr uniforme integr}), there exists a sequence
$(Q_{n_i})_{i\in\mathbb{N}}$ in $\Omega(\beta)\subset\Omega$ that
converges in $\tau$-topology to some $Q^*$ in $\mathcal{M}(P)$
satisfying
\begin{equation}\label{inequality by lsc 1}
\phi(Q^{\ast},P)\leq\lim_{i\to
+\infty}\phi(Q_{n_i},P)=\phi(\Omega,P)<\infty.
\end{equation}
It remains to prove that $Q^*$ belongs to $\Omega$. At first, we
prove that $Q^{\ast}$ belongs to $\mathcal{M}_{\mathcal{F}}$. So,
let $f$ in $\mathcal{F}$. Denote by $Q_+^*$   the nonnegative
variation and by $Q_-^*$ the nonpositive variation of $Q^*$:
$Q^*=Q_+^*-Q_-^*$. Using Fenchel's inequality through the integral
we can write
\begin{eqnarray}\label{int f d Q star 1-1}
\int |f|~dQ^*_+ & = &\int |f|q^*_+~dP \nonumber\\
& \leq & \int \varphi\left(q^*_+\right)~dP+\int \varphi^*\left(|f|\right)~dP\nonumber \\
~ &\leq & \int \varphi\left(q^*\right)~dP+\int \varphi^*\left(|f|\right)~dP \nonumber\\
& = & \phi(Q^*,P)+ \int \varphi^*\left(|f|\right)~dP,
\end{eqnarray}
and similarly
\begin{equation}\label{int f d Q star 1-2}
\int |f|~dQ^*_- \leq  \phi(Q^*,P)+ \int
\varphi^*\left(|f|\right)~dP.
\end{equation}
Hence, from (\ref{inequality by lsc 1}), (\ref{int f d Q star
1-1}) and (\ref{int f d Q star 1-2}), we deduce $\int
|f|~d|Q^{\ast}|<\infty$ since $$\int |f|~d|Q^*|=\int |f|~dQ^*_+ +
\int |f|~dQ^*_-.$$ Hence $Q^*$ belongs to
$\mathcal{M}_\mathcal{F}$. We still have to prove that $Q^*$
belongs to $\Omega$. Since $\Omega$ is, by assumption, a closed
set in $\left[\mathcal{M}_\mathcal{F};\tau_\mathcal{F}\right]$, it
is enough to show that the sequence ${\left(Q_{n_i}\right)}_{i}$
(which belongs to $\Omega(\beta)\subset \Omega$) converges to
$Q^{\ast}$ in
$\left[\mathcal{M}_{\mathcal{F}};\tau_{\mathcal{F}}\right]$. Note
that the sequence $Q_{n_i}$ converges to $Q^*$ in $\tau$-topology.
Hence, we still have to prove that $\int f~dQ_{n_i}$ converges to
$\int f~dQ^*$ for all $f$ in $\mathcal{F}$.  So, let $f$ in
$\mathcal{F}$. We use now similar argument as in the proof of
Lemma 2.1 in \cite{EischelsbacherScmock1997}. Let $\epsilon >0$.
Define $\alpha=(m+\beta)/\epsilon$. Using the fact that
$\varphi^*(0)=0$, by condition (\ref{condition Fenchel ineq}) and
the dominated convergence theorem, there exists $j_0\in\mathbb{N}$
such that
\begin{equation*}
\frac{1}{\alpha}\int \varphi^*\left(\alpha
|f|\mathds{1}_{\left\{|f|>j\right\}}\right)~dP<\epsilon
\end{equation*}
for all $j\geq j_0$. Hence, using Fenchel's inequality and the
fact that the sequence $(Q_{n_i})_i$ belongs to $\Omega(\beta)$,
we can write
\begin{eqnarray}
\left|\int f~dQ_{n_i}-\int f\mathds{1}_{\left\{|f|\leq
j\right\}}~dQ_{n_i}\right| & \leq & \int
\left|f-f\mathds{1}_{\left\{|f|\leq
j\right\}}\right|~d|Q_{n_i}|\nonumber\\
& = &  \frac{1}{\alpha}\int \alpha
|f|\mathds{1}_{\left\{|f|>j\right\}}~d|Q_{n_i}|\nonumber\\
& \leq & 2\left[\frac{1}{\alpha}\phi(Q_{n_i},P)+
\frac{1}{\alpha}\int\varphi^*\left(\alpha
|f|\mathds{1}_{\left\{|f|>j\right\}}\right)~dP\right]\nonumber\\
& \leq &
2\left[\frac{1}{\alpha}(m+\beta)+\epsilon\right]=4\epsilon.
\end{eqnarray}
We have just proved that, for all $f\in\mathcal{F}$, for all
$\epsilon >0$, there exists $j_0\in\mathbb{N}$, such that for all
$j\geq j_0$ and all $i\in\mathbb{N}$,
\begin{equation}\label{encadrement}
\int f\mathds{1}_{\left\{|f|\leq j\right\}}~dQ_{n_i} -4\epsilon
\leq \int f~dQ_{n_i} \leq \int f\mathds{1}_{\left\{|f|\leq
j\right\}}~dQ_{n_i} + 4\epsilon.
\end{equation}
Using the fact that the sequence $(Q_{n_i})_i$ converges to $Q^*$
in $\tau$-topology, by passage to limits in (\ref{encadrement})
when $i \to \infty$, then when $j\to \infty$ and finally when
$\epsilon\to 0$, we get $\lim_{i\to\infty}\int f~dQ_{n_i}=\int
f~dQ^*$. Hence, the sequence $(Q_{n_i})_i$ converges to $Q^*$ in
$\tau_\mathcal{F}$-topology, which implies that $Q^*$ belongs to
$\Omega$ since $\Omega$ is closed in
$\left[\mathcal{M}_\mathcal{F};\tau_\mathcal{F}\right]$. From the
inequality (\ref{inequality by lsc 1}), we conclude that $Q^*$ is
a $\phi$-projection of $P$ on $\Omega$. This completes  the proof.
$\blacksquare$\\

\noindent Using Hölder inequality, we give in the following
Theorem another result of existence of $\phi$-projection on closed
set in $\left[\mathcal{M}_\mathcal{F};\tau_\mathcal{F}\right]$. In
the sequel, $\left\|.\right\|_k$ denotes the usual norm of the
vector space $L_k(\mathcal{X},P)$, $1\leq k\leq +\infty$.\\

\begin{theorem}\label{thm exist proj sur ferme}
Let $\Omega$ be some closed set in $\mathcal{M}_{\mathcal{F}}$
equipped with the $\tau _{\mathcal{F}}$-topology. Assume that the
following conditions
\begin{equation}\label{phi Omega P est fini}
\phi(\Omega,P) <\infty,
\end{equation}
\begin{eqnarray}\label{condition pr existence}
\text{there exists numbers } 1< r, k < +\infty \text{ such that }
r^{-1}+k^{-1}=1,\\ \lim_{\left\vert x\right\vert \rightarrow
\infty }\frac{\varphi(x)}{\left\vert x\right\vert ^{r}}>0, ~\text{
and for every } f \in \mathcal{F}, ~ \left\|f\right\|_k<\infty
\nonumber
\end{eqnarray}
hold. Then there exists at least one $\phi$-projection of $P$ on
$\Omega$.\\
\end{theorem}

\noindent \textbf{Proof of Theorem \ref{thm exist proj sur ferme}}
Since condition (\ref{condition pr existence}) implies
(\ref{condition pr uniforme integr}), as in the proof of Theorem
\ref{corollary thm exist proj sur ferme},  there exists a sequence
$(Q_{n_i})_{i\in\mathbb{N}}$ in $\Omega(\beta)\subset\Omega$ that
converges in $\tau$-topology to some $Q^*$ in $\mathcal{M}(P)$
satisfying
\begin{equation}\label{inequality by lsc}
 \phi(Q^{\ast},P)\leq\lim_{i\to
+\infty}\phi(Q_{n_i},P)=\phi(\Omega,P)<\infty.
\end{equation}
We have to prove that $Q^*$ belongs to $\Omega$.
 At first, we prove that $Q^{\ast}$ belongs to
$\mathcal{M}_{\mathcal{F}}$. For all $f$ in $\mathcal{F}$, we have
\begin{eqnarray}\label{int f d Q star}
\int |f|~d|Q^{\ast}| & = &\int |f||q^{\ast}|~dP\nonumber\\
 &  = &
\int |f||q^{\ast}|\mathds{1}_{\left\{ |q^{\ast}|\leq c_0\right\}
}~dP+\int |f||q^{\ast}|\mathds{1}_{\left\{ |q^{\ast}|>c_0\right\}}~dP\nonumber\\
& \leq & c_0\int |f|~dP+\int |f|\frac{|q^{\ast}|}{{\varphi
(q^{\ast})}^{1/r}}{\varphi (q^{\ast})}^{1/r}
\mathds{1}_{\left\{ |q^{\ast}|>c_0\right\}}~dP \nonumber \\
~ &\leq & c_0\int |f|~dP+c_1 {\left(\int
{|f|}^{k}~dP\right)}^{1/k}{\left(\int\varphi(q^{\ast})~dP\right)}^{1/r}\nonumber\\
& = & c_0\int |f|~dP+c_1 {\left(\int
{|f|}^{k}~dP\right)}^{1/k}{\left(\phi(Q^*,P)\right)}^{1/r}.
\end{eqnarray}
Hence, from (\ref{inequality by lsc}) and (\ref{condition pr
existence}), we deduce $\int |f|~d|Q^{\ast}|<\infty$. We still
have to prove that $Q^{\ast}$ belongs to $\Omega$. Since $\Omega$
is, by assumption, a closed set in
$\left[\mathcal{M}_\mathcal{F};\tau_\mathcal{F}\right]$, it is
enough to show that the sequence ${\left(Q_{n_i}\right)}_{i}$
(which belongs to $\Omega(\beta)\subset \Omega$) converges to
$Q^{\ast}$ in
$\left[\mathcal{M}_{\mathcal{F}};\tau_{\mathcal{F}}\right]$. Note
that the sequence $Q_{n_i}$ converges to $Q^*$ in $\tau$-topology.
Hence, we still have to prove that $\int f~dQ_{n_i}$ converges to
$\int f~dQ^*$ for all $f$ in $\mathcal{F}$.  So, let $f$ in
$\mathcal{F}$. For all positive number $b$, using (\ref{condition
pr existence}), we can write $\int f~dQ_{n_i}  =  \int
f\mathds{1}_{\left\{ |f|\leq b\right\} }~dQ_{n_i}+\int
f\mathds{1}_{\left\{|f|>b\right\} }~dQ_{n_i}=:A+B,$ and
\begin{eqnarray*}
|B| & = & \left\vert\int f\mathds{1}_{\left\{ |f|>b\right\}
}~dQ_{n_i}\right\vert \leq \int |f|\mathds{1}_{\left\{
|f|>b\right\}}~d\left\vert Q_{n_i}\right\vert ~ = ~
\int |f|\mathds{1}_{\left\{ |f|>b\right\}}|q_{n_i}|~dP\\
& = & \int |f|\mathds{1}_{\left\{ |f|>b\right\}
}|q_{n_i}|\mathds{1}_{\left\{ |q_{n_i}|\leq c_0\right\} }~dP+\int
|f|\mathds{1}_{\left\{
|f|>b\right\}}|q_{n_i}|\mathds{1}_{\left\{|q_{n_i}|>c_0\right\}}~dP \\
 &\leq & c_0\int |f|\mathds{1}_{\left\{|f|>b\right\}}~dP+\int
|f|\mathds{1}_{\left\{|f|>b\right\}}\frac{|q_{n_i}|}{{\varphi
(q_{n_i})}^{1/r}}{
\varphi (q_{n_i})}^{1/r}\mathds{1}_{\left\{|q_{n_i}|>c_0\right\}}~dP \\
 & \leq & c_0\int |f|\mathds{1}_{\left\{|f|>b\right\}}~dP+c_1
\int |f|\mathds{1}_{\left\{|f|>b\right\} }{\varphi(q_{n_i})}^{1/r}~dP \\
 & \leq & c_0\int |f|\mathds{1}_{\left\{|f|>b\right\}}~dP+c_1
{\left(\int {|f|}^{k}\mathds{1}_{\left\{|f|>b\right\}}~dP\right)
}^{1/k}{\left(\int \varphi(q_{n_i})~dP\right)}^{1/r}.
\end{eqnarray*}
We deduce
\begin{equation}\label{B1 and B2}
(B1) \leq \int f~dQ_{n_i}\leq (B2),
\end{equation}
with
\begin{eqnarray}
(B1) &  :=  &  \int f\mathds{1}_{\left\{|f|\leq b\right\}
}~dQ_{n_i}-c_0\int |f| \mathds{1}_{\left\{|f|>b\right\}}~dP  -
  c_1 {\left( \int {|f|}^{k}\mathds{1}_{\left\{
|f|>b\right\} }~dP\right) }^{1/k}{\left( \int \varphi
(q_{n_i})~dP\right) }^{1/r},\nonumber
\end{eqnarray}
and
\begin{eqnarray}
(B2)  & :=  & \int f\mathds{1}_{\left\{|f|\leq b\right\}
}~dQ_{n_i}+c_0\int |f| \mathds{1}_{\left\{|f|>b\right\}}~dP +
  c_1 {\left( \int {|f|}^{k}\mathds{1}_{\left\{
|f|>b\right\}}~dP\right)}^{1/k}{\left(\int \varphi
(q_{n_i})~dP\right)}^{1/r}.\nonumber
\end{eqnarray}
The functions $\left\{f_{b}:=|f|\mathds{1}_{\left\{
|f|>b\right\}},~~b\geq 0\right\}$ and
$\left\{f_{b}^k:=|f|^k\mathds{1}_{\left\{ |f|>b\right\}},~~b\geq
0\right\}$ are dominated respectively by $|f|$ and $|f|^k$.
Moreover, $\int |f|~dP$ and $\int |f|^k~dP$ are finite by
assumption (\ref{condition pr existence}). We thus get by the
dominated convergence theorem
\begin{equation*}
\lim_{b\rightarrow +\infty}\int |f|\mathds{1}_{\left\{
|f|>b\right\}}~dP=\lim_{b\rightarrow +\infty}\int
|f|^k\mathds{1}_{\left\{ |f|>b\right\}}~dP=0.
\end{equation*}
Hence, from (\ref{B1 and B2}), we get
\begin{equation*}
\int f~dQ^{\ast}=\lim_{b\rightarrow +\infty}\lim_{i\rightarrow
+\infty}(B1)\leq \lim_{i\rightarrow +\infty}\int f~dQ_{n_i}\leq
\lim_{b\rightarrow +\infty}\lim_{i\rightarrow +\infty}(B1)=\int
f~dQ^{\ast},
\end{equation*}
which is to say that the subsequence ${\left( Q_{n_i}\right)}_{i}$
converges to $Q^{\ast}$ in $\tau_{\mathcal{F}}$-topology. Hence,
$Q^{\ast}$ belongs to $\Omega$. From inequality (\ref{inequality
by lsc}), we deduce that $Q^*$ is a $\phi$-projection of $P$ on
$\Omega$. This ends the proof of Theorem
\ref{thm exist proj sur ferme}. $\blacksquare$\\

\noindent Note that the above results do not apply in the case of
$KL_m$ and Hellinger divergences since the condition
$\lim_{|x|\to\infty}\frac{\varphi(x)}{|x|}=+\infty$ does not hold.
The following Corollary applies without assumption
$\lim_{|x|\to\infty}\frac{\varphi(x)}{|x|}=+\infty$, in
particular, in the $KL_m$ and Hellinger divergences cases.\\

\begin{corollary}\label{existence de proj pour klm and h}
Let $\Omega$ be a closed set in $\left[\mathcal{M},\tau\right]$.
If the following condition:  there exists
\begin{equation}\label{une condition pr unifor integ}
 u,l\in L_1(\mathcal{X},P) \text{ such that }
 u\leq \frac{dQ}{dP}\leq l ~(P-a.e.) \text{ for all }
 Q\in\Omega\cap\text{dom}\phi
\end{equation}
holds, then there exists at least one $\phi$-projection of $P$ on
$\Omega$ whenever $\phi(\Omega,P)$ is finite.\\
\end{corollary}

\noindent \textbf{Proof of Corollary \ref{existence de proj pour
klm and h}} Similar to that of Theorem \ref{corollary thm exist
proj sur ferme}. The uniform integrability of the set
$\Lambda(\beta)$ holds by condition (\ref{une condition pr unifor
integ}). $\blacksquare$

\section{Characterization of
$\phi$-Projections on sets of signed measures} \noindent In this
section, we extend known results pertaining to the
characterization of the $\phi$-projections as can be found in
\cite{Ruschendorf1984}, \cite{Ruschendorf1987},
\cite{Liese-Vajda1987}, (see also \cite{Csiszar1975} and
\cite{Csiszar1984} for the characterization of $KL$-projections).
These authors have characterized the $\phi$-projections on subsets
of $\mathcal{M}^{1}$. We expose similar results when considering
subsets of $\mathcal{M}$ and take the occasion to clarify some
proofs. We first consider the case of general subsets $\Omega$ of
$\mathcal{M}$ and then the case of convex subsets of $\mathcal{M}$
defined by linear constraints. For the whole Section, we assume
that the convex function $\varphi$ is differentiable.
\subsection{On general Sets $\Omega$}
\noindent We will use the following assumption
\begin{equation}\label{condition C.0}
\begin{tabular}{l}
 There exists  $0<\delta<1$ such that for all $c$ in
$\left[1-\delta ,1+\delta \right]$, \\
we can find  numbers $c_{1},c_{2},$ $c_{3}$ such that \\
$\varphi (cx)\leq c_{1}\varphi (x)+c_{2}\left\vert x\right\vert
+c_{3}$, for all real $x$.
\end{tabular}
\end{equation}

\begin{remark}
Condition (\ref{condition C.0}) holds for all power divergences
including $KL$, $KL_m$ and Hellinger divergences. Note also that
condition (\ref{condition C.0}) implies that $a_\varphi$ equals
$0$ or $-\infty$ and $b_\varphi$ equals $+\infty$.\\
\end{remark}

\begin{remark}
In all the sequel, condition (\ref{condition C.0}) above can be
replaced by any other condition which implies part (1) of Lemma
\ref{lemma caract proj 1} below.\\
\end{remark}

\noindent We first give two Lemmas, which we will use in the proof
of Theorem \ref{thm caract proj} and Theorem \ref{thm rusch}
below.\\

\begin{lemma}\label{lemma caract proj 1}
Assume that (\ref{condition C.0}) holds. Then, for all $Q$ in
$\mathcal{M}$ such that $\phi(Q,P)$ is finite, we have
\begin{enumerate}
 \item  [(1)] for any $c$ in $[1-\delta ,1+\delta ]$, $\varphi\left(
c\frac{dQ}{dP}\right)$ belongs to $L_1(\mathcal{X},P)$.
 \item  [(2)] $\lim_{c\uparrow 1}\phi (cQ,P)=\phi (Q,P)=\lim_{c\downarrow
1}\phi(cQ,P)$.\\
\end{enumerate}
\end{lemma}

\noindent \textbf{Proof of Lemma \ref{lemma caract proj 1}} (1)
Under condition (\ref{condition C.0}), for all  $Q$ in
$\mathcal{M}$ such that $\phi(Q,P)<\infty$, we have
\begin{equation*}
\varphi \left(c\frac{dQ}{dP}\right) \leq c_{1}\varphi \left(
\frac{dQ}{dP}\right) +c_{2}\left\vert\frac{dQ}{dP}\right\vert
+c_{3}.
\end{equation*}
Integrating with respect to $P$ yields
\begin{equation*}
\int \varphi\left(c\frac{dQ}{dP}\right)~dP\leq c_{1}\phi
(Q,P)+c_{2}\int \left\vert \frac{dQ}{dP}\right\vert~
dP+c_{3}<\infty.
\end{equation*}
(2) For all $c$ in $[1-\delta ,1+\delta ]$, define the functions
\begin{equation*}
\begin{array}{ccccccc}
l_{c} & : & x\in\mathbb{R} & \mapsto & l_{c}(x) & := &
\varphi (cx)\mathds{1}_{]-\infty ,0[}(cx), \\
g_{c} & : & x\in\mathbb{R} & \mapsto & g_{c}(x) & := &
\varphi (cx)\mathds{1}_{[0,1]}(cx), \\
h_{c} & : & x\in\mathbb{R} & \mapsto & h_{c}(x) & := & \varphi
(cx)\mathds{1}_{]1,+\infty[}(cx).
\end{array}
\end{equation*}
For any $c$ and $x$, we have $\varphi
(cx)=l_{c}(x)+g_{c}(x)+h_{c}(x)$. For all real $x$, the functions
$c\rightarrow l_{c}(x)$ and $c\rightarrow h_{c}(x)$ are
 nondecreasing, and the function $c\rightarrow g_{c}(x)$ is nonincreasing.
Denote $q:=\frac{dQ}{dP}$. Apply the monotone convergence theorem
to get
\begin{equation*}
\lim_{c\uparrow 1}\int l_{c}(q)~dP=\int l_{1}(q)~dP ~\text{ and }~
\lim_{c\uparrow 1}\int h_{c}(q)~dP=\int h_{1}(q)~dP.
\end{equation*}
On the other hand, the class of functions $\left\{x\rightarrow
g_{c}(x), ~c\text{ in } [1-\delta ,1+\delta ]\right\}$ is bounded
above by the function $x\rightarrow g_{1-\delta }(x)$.
Furthermore, for all $Q$ in $\mathcal{M}$, $g_{1-\delta}(q)$
belongs to $L_1(\mathcal{X},P)$ by the condition (\ref{condition
C.0}). Hence, applying the monotone convergence theorem we get
\begin{equation*}
\lim_{c\uparrow 1}\int g_{c}(q)~dP=\int g_{1}(q)~dP.
\end{equation*}
Those three limits prove the first part of the claim. The same
argument completes the proof of the Lemma. $\blacksquare$\\

\begin{lemma}\label{lemma caract proj}
Assume that condition (\ref{condition C.0}) holds. Then, for all
$Q$ in $\text{dom}\phi$, $\varphi^{\prime}(q)q$ belongs to
$L_1(\mathcal{X},P)$, where $q:=\frac{dQ}{dP}$.\\
\end{lemma}

\noindent\textbf{Proof of Lemma \ref{lemma caract proj}} Using the
convexity of the function $\varphi$, for all
 $\epsilon > 0$, we have
\begin{equation*}
\frac{\varphi(q)-\varphi \left((1-\epsilon)q\right)}{\epsilon}\leq
q\varphi^{\prime}(q)\leq\frac{\varphi\left((1+\epsilon)q\right)
-\varphi(q)}{\epsilon}.
\end{equation*}
By Lemma \ref{lemma caract proj 1}, for all $\epsilon$ satisfying
$0 < \epsilon<\delta$, both the LHS and the RHS terms belong to
$L_1(\mathcal{X},P)$, and hence $\varphi^{\prime}(q)q\in
L_1(\mathcal{X},P)$. $\blacksquare$\\

\begin{theorem}\label{thm caract proj}
Let $\Omega$ be a subset of $\mathcal{M}$ and $Q^*$ be a signed
measure in $\Omega\cap\text{dom}\phi$. Then
\begin{enumerate}
 \item [(1)] The following are sufficient conditions for $Q^*$ to be a $\phi$-projection of
 $P$ on $\Omega$:
  (i) $\varphi^{\prime}(q^{\ast})q\in L_1(\mathcal{X},P)$ and (ii)
$\int\varphi^{\prime}(q^{\ast})~dQ^{\ast}\leq \int
\varphi^{\prime}(q^{\ast})~dQ$, for all $Q$ in
$\Omega\cap\text{dom}\phi$.
 \item [(2)] If condition (\ref{condition C.0}) holds and $\Omega$ is
 convex, then these conditions are necessary as well.\\
\end{enumerate}
\end{theorem}

\noindent \textbf{Proof of Theorem \ref{thm caract proj}}
Convexity and differentiability of $\varphi$ imply, for all
positive $\epsilon$,
\begin{equation}\label{ineg caract proj}
\varphi^{\prime}(q^{\ast})(q-q^{\ast})\leq \frac{\varphi\left(
(1-\epsilon)q^{\ast}+\epsilon q\right)
-\varphi(q^{\ast})}{\epsilon} \leq \varphi(q)-\varphi(q^{\ast}).
\end{equation}
The middle term in the above display, by the convexity of
$\varphi$, decreases to $\varphi^{\prime}(q^{\ast})(q-q^{\ast})$
when $\epsilon \downarrow 0$. Furthermore, it is bounded above by
$\varphi(q)-\varphi(q^{\ast})$ which belongs to
$L_1(\mathcal{X},P)$ for all $Q$ in $\text{dom}\phi$. Hence,
applying the monotone convergence theorem to get
\begin{equation}\label{a m c t to get}
\int \varphi^{\prime}(q^{\ast})(q-q^{\ast})~dP=\lim_{\epsilon
\downarrow 0} \int
\frac{\varphi\left((1-\epsilon)q^{\ast}+\epsilon q\right)
-\varphi(q^{\ast})}{\epsilon}~dP, ~~\text{ for all
}~Q\in\text{dom}\phi.
\end{equation}
Proof of part (1): Integrating $\left(\ref{ineg caract
proj}\right)$ with respect to $P$ and using (i) and (ii) in part
(1) of the Theorem, we obtain for all $Q$ in
$\Omega\cap\text{dom}\phi$
\begin{equation}
\phi(Q,P)-\phi(Q^{\ast},P)\geq
\int\varphi^{\prime}(q^{\ast})(q-q^{\ast})~dP=\int \varphi
^{\prime}(q^{\ast})~dQ-\int\varphi^{\prime}(q^{\ast})~dQ^{\ast}\geq
0.  \label{ineg phi prime}
\end{equation}
Hence, $Q^{\ast}$ is a $\phi$-projection of $P$ on $\Omega$. Proof
of part (2):  Convexity of both $\Omega$ and $\text{dom}\phi$,
implies that for all $Q\in\Omega\cap\text{dom}\phi$,
$(1-\epsilon)Q+\epsilon Q^*$ belongs to
$\Omega\cap\text{dom}\phi$. Since $Q^*$ is a $\phi$-projection of
$P$ on $\Omega$,  for all $Q\in\Omega\cap\text{dom}\phi$ and all
$\epsilon$ satisfying $0<\epsilon<1$, we get
    $\phi\left((1-\epsilon)Q+\epsilon Q^*, P\right)-\phi(Q^*,
    P)\geq 0.$
Combining this with (\ref{a m c t to get}) and using the fact that
$Q^*$ is a $\phi$-projection of $P$ on $\Omega$, we obtain for all
$Q$ in $\Omega\cap\text{dom}\phi$
\begin{eqnarray} \label{ineg 1}
\int \varphi^{\prime}(q^{\ast})(q-q^{\ast})~dP & = &
\lim_{\epsilon \downarrow 0} \int
\frac{\varphi\left((1-\epsilon)q^{\ast}+\epsilon q\right)
-\varphi(q^{\ast})}{\epsilon}~dP \nonumber\\
 & = & \lim_{\epsilon \downarrow 0}\frac{1}{\epsilon}\left[
 \phi\left((1-\epsilon)Q^*+\epsilon Q, P\right)-\phi\left(Q^*,
 P\right)\right] \geq 0.
\end{eqnarray}
On the other hand, integrating (\ref{ineg caract proj}) with
respect to $P$, we obtain for all $Q$ in
$\Omega\cap\text{dom}\phi$
\begin{equation}\label{ineg 2}
   \int \varphi^{\prime}(q^{\ast})(q-q^{\ast})~dP \leq
   \phi(Q,P)-\phi(Q^*,P)<\infty.
\end{equation}
Hence, (\ref{ineg 1}) and (\ref{ineg 2}) imply
\begin{equation}\label{phiprime q satr q-q star}
\varphi^{\prime}(q^{\ast})(q-q^{\ast})\in L_1(\mathcal{X},P),
~~~\text{ for all }~Q\in \Omega\cap\text{dom}\phi.
\end{equation}
By Lemma \ref{lemma caract proj},
$\varphi^{\prime}(q^{\ast})q^{\ast}\in L_1(\mathcal{X},P)$.
Combining this with (\ref{phiprime q satr q-q star}), we obtain
that
\begin{equation*}
    \text{for all }~Q\in \Omega\cap\text{dom}\phi,~\text{ we have
    }~\varphi'(q^*)q\in L_1(\mathcal{X},P)
\end{equation*}
and $\int    \varphi'(q^*)~dQ^* \leq \int \varphi'(q^*)~dQ.$ This
completes the proof of  Theorem \ref{thm caract proj}.
$\blacksquare$

\subsection{On Sets defined by Linear Constraints}
\noindent In this subsection, we consider the problems of
existence and characterization of $\phi$-projections of some p.m.
$P$ on linear set $S$ of measures in $\mathcal{M}$ defined by
arbitrary family of constraints. So, let $\mathcal{G}$ denote a
collection (finite or infinite, countable or not) of real valued
functions defined on $\left(\mathcal{X},\mathcal{B}\right)$. The
class $\mathcal{G}$ is assumed to contain the function
$\mathds{1}_\mathcal{X}$. The set $S$ is defined by
\begin{equation}\label{Omega def par contr lin f s m}
S :=\left\{ Q\in \mathcal{M}_\mathcal{G}(P)   \text{ such that }
\int_\mathcal{X} dQ=1, \int_\mathcal{X} g~dQ=0, \text{ for all }
g\text{ in }\mathcal{G\setminus
}\left\{{\mathds{1}}_\mathcal{X}\right\}\right\}.
\end{equation}
The following result states the explicit form of $Q^{\ast}$, a
$\phi$-projection of $P$ on $S$, when it exists.\\

\begin{theorem}\label{thm rusch}\
\begin{enumerate}
\item [(1)] Let $Q^*$ be some finite measure in
$S\cap\text{dom}\phi$. A sufficient condition, for $Q^*$ to be a
$\phi$-projection of $P$ on $S$,  is that  there exists numbers
$c_1,\ldots,c_d\in\mathbb{R}$ and functions
$g_1,\ldots,g_d\in\mathcal{G}$ such that $\varphi^{\prime
}(q^{\ast}(x))=c_1g_1(x)+\cdots+c_dg_d(x)$~ ($P$-a.e.).
 \item [(2)] Assume that condition (\ref{condition C.0}) holds. Then, any
 $\phi$-projection, say $Q^*$, of $P$ on $S$,  if it exists,
 satisfies
$\varphi^{\prime}(q^{\ast})$ belongs to $\overline{\left\langle
\mathcal{G}\right\rangle }$, (the closure of $\left\langle
\mathcal{G}\right\rangle$) in $L_1(\mathcal{X},|Q^{\ast}|)$.\\
\end{enumerate}
\end{theorem}

\noindent If $\mathcal{G}$ is a finite collection  of functions in
$L_1(\mathcal{X},|Q^*|)$, then the vector space $\langle
\mathcal{G}\rangle$ is closed in $L_1(\mathcal{X},|Q^*|)$. So,
from the above Theorem, we can state the following Corollary:\\

\begin{corollary}\label{finite number of constr - corollary1}
 Let $\mathcal{G}:=\left\{\mathds{1}_\mathcal{X},g_1,
 \ldots,g_l\right\}$
 be a finite collection of measurable functions on $\mathcal{X}$. Then
 (1) and (2) below hold.
\begin{enumerate}
 \item [(1)] Let $Q^*$ be some measure in $S\cap\text{dom}\phi$. A sufficient condition, for $Q^*$ to be
 a $\phi$-projection of $P$ on $S$, is that there exists some
 constant $c\in \mathbb{R}^{1+l}$ such that
 \begin{equation*}
\varphi'\left(\frac{dQ^*}{dP}(x)\right)=c_0+\sum_{i=1}^{l}c_ig_i(x)\quad
(P-a.e.).
 \end{equation*}
  \item  [(2)] Assume that condition (\ref{condition C.0}) holds. Then any $\phi$-projection, say $Q^*$, of $P$ on $S$, if
  it exists, satisfies
  \begin{eqnarray}
  \text{there exists some constant } c\in \mathbb{R}^{1+l} \text{ such that
  }\nonumber\\
  \varphi'\left(\frac{dQ^*}{dP}(x)\right)=c_0+\sum_{i=1}^{l}c_ig_i(x)\quad
(|Q^*|-a.e.).\nonumber
  \end{eqnarray}\\
\end{enumerate}
\end{corollary}

\noindent It should be noticed that the preceding Theorem and
Corollary do not provide a definite description of the projected
measure; indeed, it does not give any information on the support
of $|Q^\ast|$ (see example \ref{contre exemple} below). However,
if $\varphi(0)=+\infty$ (which holds for example for the
$KL_m$-divergence), then any $\phi$-projection $Q^*$ of $P$ on
some set $\Omega$, if it exists, has obviously the same support as
$P$ when $\phi(\Omega,P)$ is finite. Furthermore, we prove in the
following Lemma that if $\varphi'(0)=-\infty$ (which holds for
instance in the case of $KL$, $KL_m$ and Hellinger divergences),
then any $\phi$-projection of $P$ on some convex set $\Omega$ when
it exists has the same support as $P$. At first, state the
following Corollary which applies in the $KL_m$-divergence case.\\

\begin{corollary} Let $\mathcal{G}$ be defined as in
Corollary \ref{finite number of constr - corollary1}. Assume that
assumption (\ref{condition C.0}) holds. Suppose additionally that
$\varphi(0)=+\infty$, and let $Q^*$ be some p.m. in
$S\cap\text{dom}\phi$. Then
 $Q^*$ is a $\phi$-projection of $P$ on
$S$ iff there exists some constant $c\in \mathbb{R}^{1+l}$ such
that
 \begin{equation*}
\varphi'\left(\frac{dQ^*}{dP}(x)\right)=c_0+\sum_{i=1}^{l}c_ig_i(x)\quad
(P-a.e.).
\end{equation*}
\end{corollary}

\begin{lemma}\label{meme support}
Assume that condition (\ref{condition C.0}) holds, $a_\varphi=0$
and  $\varphi'(0)=-\infty$. Let $\Omega$ be some convex set of
signed finite measures. If there exists some
$Q_0\in\Omega\cap\text{dom}\phi$ such that $\frac{dQ_0}{dP}>0$~
($P$-a.e.), then any $\phi$-projection, say $Q^*$, of $P$ on
$\Omega$, if it exists, has the same support as $P$, i.e.,
$\frac{dQ^*}{dP}>0$ ~ ($P$-a.e.).\\
\end{lemma}

\noindent \textbf{Proof of Lemma \ref{meme support}} Let
$A:=\left\{x\in\mathcal{X}; ~q^*(x)=0\right\}$. Suppose that
$P(A)>0$. Since $Q_0$ and $P$ have the same support by assumption,
$Q_0(A)>0$.  By (\ref{ineg caract proj}) (replacing $Q$ by $Q_0$),
$Q_0(A)>0$ implies that $\int \varphi'(q^*)q~dP=-\infty$ since
$\int \left|\varphi'(q^*)q^*\right|~dP<\infty$. This contradicts
(\ref{ineg 1}),
which completes the proof.$\blacksquare$\\
\\
\noindent We can  now state, from the above Theorem, the following
Corollary which applies in the case of $KL$, $KL_m$ and Hellinger
divergences.\\

\begin{corollary}
Let $\mathcal{G}$ be defined as in Corollary \ref{finite number of
constr - corollary1}. Assume that assumption (\ref{condition C.0})
holds. Suppose additionally that $a_\varphi=0$ and
$\varphi'(0)=-\infty$. If there exists some $Q_0\in
S\cap\text{dom}\phi$ such that $\frac{dQ_0}{dP}>0$~ ($P$-a.e.),
then the following holds : a p.m. $Q^*$ in $S\cap\text{dom}\phi$
is a $\phi$-projection of $P$ on $S$ iff there exists some
constant $c\in \mathbb{R}^{1+l}$ such that
 \begin{equation*}
\varphi'\left(\frac{dQ^*}{dP}(x)\right)=c_0+\sum_{i=1}^{l}c_ig_i(x)\quad
(P-a.e.).
\end{equation*}\\
\end{corollary}

\begin{remark}
Versions of Theorem \ref{thm rusch}, for sets of p.m.'s, have been
proved   by \cite{Csiszar1975} and \cite{Csiszar1984} for the
Kullback-Leibler divergence, and by \cite{Ruschendorf1984} and
\cite{Liese-Vajda1987} for $\phi$-divergences between p.m.'s. We
prove it in the present context, that is when the set $S$ (see
(\ref{Omega def par contr lin f s m})) is a subset of signed
finite measures and $P$ is a p.m. using similar techniques.\\
\end{remark}

\noindent \textbf{Proof of Theorem \ref{thm rusch}}   We start by
proving (1). If $\varphi^{\prime}(q^{\ast})$ belongs to
$\langle\mathcal{G}\rangle$, then for all $Q$ in
 $S$, we have  $\int\varphi^{\prime}(q^{\ast})~dQ^{\ast}=\int\varphi^{\prime}(q^{\ast})~dQ$
which, by the first part of Theorem \ref{thm caract proj}, proves
that $Q^{\ast}$ is a $\phi$-projection of $P$ on $ S$. Proof of
part (2): Since $Q^{\ast}$ is a signed finite measure, by the Hahn
decomposition theorem, there exists a partition
$\mathcal{X}=\mathcal{X}_1\cup\mathcal{X}_2$ such that
$\mathcal{X}_1,\mathcal{X}_2\in \mathcal{B}$ and satisfying
\begin{enumerate}
 \item [] for all $B\in\mathcal{B}$,  such that
$B\subset\mathcal{X}_1$ we have $Q^\ast (B)\geq 0$
\end{enumerate}
and
\begin{enumerate}
 \item [] for all $B\in\mathcal{B}$,  such that
$B\subset\mathcal{X}_2$ we have $Q^\ast (B)\leq 0$.
\end{enumerate}
 Denote by
$Q^\ast_+$ and $Q^\ast_-$ respectively the nonnegative variation
and  the nonpositive variation of $Q^\ast$ which are defined, for
all $B\in\mathcal{B}$, by
\begin{equation*}
Q^\ast_+(B):=Q^\ast(B\cap\mathcal{X}_1) ~\text{ and
}~Q^\ast_-(B):=-Q^\ast(B\cap\mathcal{X}_2).
\end{equation*}
So, $Q^\ast_+$ and $Q^\ast_-$ are nonnegative finite measures,
$Q=Q^\ast_+-Q^\ast_-$ and the total variation $|Q^\ast|$ is, by
definition, the nonnegative measure $Q^\ast_+ +Q^\ast_-$. Denote
by ${\langle\mathcal{G}\rangle}^\perp_+$ and by
${\langle\mathcal{G}\rangle}^\perp_-$  respectively the orthogonal
of $\langle\mathcal{G}\rangle$ in $L_1(\mathcal{X},Q^\ast_+)$ and
in $L_1(\mathcal{X},Q^\ast_-)$, i.e., the sets defined by
\begin{equation*}\label{orthogonal de G plus}
    {\langle\mathcal{G}\rangle}^\bot_+:=\left\{h\in L_\infty(\mathcal{X},Q^*_+)
    ~\text{ such that }~ \int f h~dQ^*_+=0, ~\text{ for all }~
    f\in\langle\mathcal{G}\rangle \right\}
\end{equation*}
and
\begin{equation*}\label{orthogonal de G moins}
    {\langle\mathcal{G}\rangle}^\bot_-:=\left\{h\in L_\infty(\mathcal{X},Q^*_-)
    ~\text{ such that }~ \int f h~dQ^*_-=0, ~\text{ for all }~
    f\in\langle\mathcal{G}\rangle \right\}.
\end{equation*}
We will prove that the two following assertions hold
\begin{equation}\label{Q star plus}
    \text{ for all  }~h\in{\langle\mathcal{G}\rangle}^\bot_+, ~\text{ we
    have }~\int\varphi'(q^*)h~dQ^*_+=0
\end{equation}
and
\begin{equation}\label{Q star moins}
    \text{ for all  }~h\in{\langle\mathcal{G}\rangle}^\bot_-, ~\text{ we
    have }~\int\varphi'(q^*)h~dQ^*_-=0.
\end{equation}
We prove (\ref{Q star plus}) by deriving a contradiction: assume
that there exists $h$ in ${\langle\mathcal{G}\rangle}^\perp_+ $
such that
  $\int\varphi^\prime(q^\ast)h~dQ^\ast_+\neq 0$.
 We then have either (a) $\int\varphi^{\prime}(q^{\ast })h~dQ^*_+<0$
or (b) $\int\varphi'(q^*)h~dQ^*_+>0$. Assume (a). For $0<\epsilon
<\delta$,\footnote{ here $\delta$ is defined in the condition
(\ref{condition C.0}).} define the measure $Q_{0}$ by
\begin{equation}\label{the measure Q zero}
dQ_{0}:=\left(1+\epsilon
\frac{h\mathds{1}_{\mathcal{X}_1}}{{\|h\|}_\infty}\right)dQ^*.
\end{equation}
Then $Q_0$ belongs to $ S$, and, following condition
(\ref{condition C.0}), $Q_0$ belongs to $\text{dom}\phi$ by Lemma
\ref{lemma caract proj 1}. Furthermore,
\begin{equation*}
\int \varphi'(q^*)~dQ_0=\int\varphi'(q^*)~dQ^*
+\epsilon\frac{1}{{\|h\|}_\infty}\int\varphi'(q^*)h~dQ^*_+<\int
\varphi'(q^*)~dQ^*,
\end{equation*}
which contradicts the fact that $Q^\ast$ is a $\phi$-projection of
$P$ on $ S$ (see part 2 in Theorem \ref{thm caract proj}). Assume
(b). Consider $-h$ instead of $h$. We thus have proved (\ref{Q
star plus}). The same arguments hold for the proof of (\ref{Q star
moins}). Therefore, $\varphi'(q^*)$ belongs to
${\left({\langle\mathcal{G}\rangle}^\bot_+\right)}^\bot_+$ and to
${\left({\langle\mathcal{G}\rangle}^\bot_-\right)}^\bot_-$
respectively the orthogonal of
${\langle\mathcal{G}\rangle}^\bot_+$ in $L_1(\mathcal{X},Q^*_+)$
and the orthogonal of ${\langle\mathcal{G}\rangle}^\bot_-$ in
$L_1(\mathcal{X},Q^*_-)$. By Hahn-Banach Theorem (see e.g. Section
2 of \cite{Brezis1983}), we have
\begin{equation*}
  {\left({\langle\mathcal{G}\rangle}^\bot_+\right)}^\bot_+=
  \overline{\langle\mathcal{G}\rangle}_+ ~\text{ and }~
   {\left({\langle\mathcal{G}\rangle}^\bot_-\right)}^\bot_-=
  \overline{\langle\mathcal{G}\rangle}_-
\end{equation*}
which are respectively the closure of $\langle\mathcal{G}\rangle$
in $L_1(\mathcal{X},Q^*_+)$ and the closure of
$\langle\mathcal{G}\rangle$ in $L_1(\mathcal{X},Q^*_-)$. This
implies that
$\varphi'(q^*)\in\overline{\langle\mathcal{G}\rangle}$ that is,
$\varphi'(q^*)$ belongs to the closure of
$\langle\mathcal{G}\rangle$ in $L_1(\mathcal{X},|Q^*|)$. This
completes  the proof of Theorem \ref{thm rusch}. $\blacksquare$\\

\begin{example}\label{contre exemple}
Let $\mathcal{X}:=[0,1]$,  $P$ be the uniform distribution on
$[0,1]$ and $\mathcal{G}:=\left\{\mathds{1}_{[0,1]},I_d\right\}$
where $I_d$ is the identity function. Consider the
$\chi_+^{2}$-divergence associated to the convex function
\begin{equation*}
\varphi(x)=\left\{
\begin{array}{lll}
\frac{1}{2}(x-1)^2 & \text{ if } & x\in[0,\infty[\\
+\infty & \text{ if } & x\in]-\infty,0[,
\end{array}
\right.
\end{equation*}
and  consider the set $M$  defined by
\begin{equation*}
M:=\left\{Q\in \mathcal{M}\text{ such that }~\int dQ=1\mbox{ and }
\int (x-1/4)~dQ(x)=0\right\}.
\end{equation*}
We apply the preceding results pertaining to the characterization
of the projection of $P$ on $M$. By Theorem \ref{thm exist proj
sur ferme}, there exists a $\chi_+^{2}$-projection, say $Q^*_+$,
of $P$ on $M$. By Theorem \ref{thm rusch}, there exists two real
numbers $c_{0}$ and $c_{1}$ such that
\begin{equation}\label{dQstar}
\frac{dQ^{\ast}_+}{dP}(x)\mathds{1}_{\left\{
q^{\ast}(x)>0\right\}}=c_{0}+c_{1}x.
\end{equation}
The support of $Q^*_+$ is different from the support of $P$; it is
strictly included in $[0,1]$. Indeed, if the support of $Q^*_+$ is
$[0,1]$, then
\begin{equation*}
dQ^\ast_+(x)=\left(c_{0}+c_{1}x\right)~dP(x)=
\left(c_{0}+c_{1}x\right)\mathds{1}_{[0,1]}(x)~dx.
\end{equation*}
Using the fact that $Q^*_+$ belongs to $M$, we obtain that
$c_0=5/2$ and $c_1=-3$. So, $Q^*_+$ satisfying
$dQ^{\ast}_+(x)=\left( 5/2-3x\right)~dP(x)$ does not belong to
$\text{dom}\chi^2_+$ (it is not a p.m.), a contradiction with the
existence of the projection. This proves that the support of
$Q^*_+$ is strictly
included in $[0,1]$.\\
\noindent Consider now the $\chi^2$-divergence, i.e., the
divergence associated to the convex function
\begin{equation*}
x\in ]-\infty,+\infty[ \mapsto \varphi(x)=\frac{1}{2}(x-1)^2,
\end{equation*}
and the set $M^1$ defined by
\begin{equation*}
M^1:=\left\{Q\in \mathcal{M}^1\text{ such that }~\int dQ=1\mbox{
and } \int (x-1/4)~dQ(x)=0\right\}.
\end{equation*}
Note that minimizing $\chi^2(.,P)$ on $M^1$ is equivalent to
minimizing $\chi^2_+(.,P)$ on $M$. Hence, $Q^*_+$ is the
$\chi^2$-projection of $P$ on $M^1$, it has not the same support
as $P$ and (\ref{dQstar}) is not a definite description of the
projection. On the other hand, the $\chi^2$-projection, say $Q^*$,
of $P$ on $M$ exists, it has the same support as $P$, it is a
signed measure  and it is characterized by $dQ^*(x)=\left(
5/2-3x\right)~dP(x)$. This example shows the interest of enhancing
$M^1$ to $M$.
\end{example}

\section{Fenchel duality for $\phi$-Divergences}
\noindent We refer to  \cite{Fenchel1949}, \cite{Moreau1962},
\cite{Brondsted1964}, \cite{Rockafellar1968},
\cite{Rockafellar1974} and \cite{EkelandTemam1999}  for the notion
of Fenchel duality of general convex functions on general vector
spaces. We consider this notion for $\phi$-divergences functionals
$Q\mapsto \phi(Q,P)$ viewed as convex functions on the vector
space of signed finite measures $\mathcal{M}_\mathcal{F}$; we give
different versions of dual representations of the
$\phi$-divergences (see Theorems \ref{duality for divergences},
\ref{representation duale v2}, \ref{representation duale v3} and
\ref{representation duale} below).
 In view of Proposition \ref{evthlc}, we identify the
topological dual space of
$\left[\mathcal{M}_\mathcal{F};\tau_{\mathcal{F}}\right]$ with
$\langle\mathcal{F}\cup\mathcal{B}_b\rangle$ and the topological
dual space of
$\left[\langle\mathcal{F}\cup\mathcal{B}_b\rangle;\tau_{\mathcal{M}}\right]$
with $\mathcal{M}_\mathcal{F}$. Hence, the Fenchel-Legendre
transform (i.e., the conjugate) of the convex function $Q\in
\left[\mathcal{M}_\mathcal{F};\tau_\mathcal{F}\right] \mapsto
\phi(Q,P)\in [0,+\infty]$ is defined as follows
\begin{equation}\label{the conjugate of phi}
f\in \left[\langle\mathcal{F}\cup\mathcal{B}_b\rangle;
\tau_\mathcal{M}\right] \mapsto
\phi^*(f):=\sup_{Q\in\mathcal{M}_\mathcal{F}}\left\{\int f
~dQ-\int \varphi\left(\frac{dQ}{dP}\right)~dP\right\},
\end{equation}
which is convex and lower semi-continuous\footnote{Note that the
conjugate of a convex function is always l.s.c. w.r.t. the weak
topology.} w.r.t. the $\tau_\mathcal{M}$-topology, the weak
topology induced on $\langle\mathcal{F}\cup\mathcal{B}_b\rangle$
by
$\mathcal{M}_\mathcal{F}$.\\
\\
\noindent By the lower semi-continuity of the convex function
$Q\in \left[\mathcal{M}_\mathcal{F};\tau_\mathcal{F}\right]
\mapsto \phi(Q,P)\in [0,+\infty]$ (see Proposition \ref{prop div
lsc param} above), applying the Fenchel duality theory (see e.g.
\cite{Rockafellar1968}, \cite{Fenchel1949} or
\cite{Dembo-Zeitouni1998} Lemma 4.5.8), we can state the following
result for any $\phi$-divergence.\\

\begin{theorem}\label{duality for divergences}
The function $Q\in
\left[\mathcal{M}_\mathcal{F};\tau_\mathcal{F}\right] \mapsto
\phi(Q,P)\in [0,+\infty]$ is the conjugate of its conjugate $f\in
\left[\langle\mathcal{F}\cup\mathcal{B}_b\rangle;
\tau_\mathcal{M}\right] \mapsto \phi^*(f)$ defined by (\ref{the
conjugate of phi}). In other words, the $\phi$-divergence
$\phi(Q,P)$ admits the dual representation
\begin{equation}\label{dual representation}
 \phi(Q,P) = \sup_{f\in
 \langle\mathcal{F}\cup\mathcal{B}_b\rangle}\left\{\int f~dQ -
 \phi^*(f)\right\}, \text{ for all }
Q\in\mathcal{M}_\mathcal{F},
\end{equation}
where $\phi^*(.)$ is defined by (\ref{the conjugate of phi}).\\
\end{theorem}

\noindent We now turn to  the calculation of $\phi^*(f)$ (in
particular the equality $\phi^*(f)=\int \varphi^*(f)~dP$), and the
problems of existence, uniqueness and characterization of a dual
optimal solution in (\ref{dual representation}) (i.e., a function
$\overline{f}\in \langle\mathcal{F}\cup\mathcal{B}_b\rangle$ which
realizes the supremum in (\ref{dual representation})) when
$\phi(Q,P)$ is
finite.\\
\\
In the following Proposition, when $\varphi$ is strictly convex
and differentiable, we give the explicit form of $\phi^*(f)$ for
all $f\in
 \langle\mathcal{F}\cup\mathcal{B}_b\rangle$ such that
 $\text{Im}f\subseteq \text{Im}\varphi'$.\\

\begin{proposition}\label{culcul de phi*}
Assume that $\varphi$ is strictly convex and differentiable, and
that for all $f,g\in\langle\mathcal{F}\cup\mathcal{B}_b\rangle$
such that $\text{Im}f\subseteq\text{Im}\varphi'$,  the integrals
\begin{equation}\label{condition pr calculer phi*}
\int |g|\left|{\varphi'}^{-1}(f)\right|~dP~\text{ and }~ \int
\varphi\left({\varphi'}^{-1}(f)\right)~dP \text{ are finite.}
\end{equation}
Then for all $f\in\langle\mathcal{F}\cup\mathcal{B}_b\rangle
\text{ such that } \text{Im}f\subseteq\text{Im}\varphi'$, we have
\begin{equation}
\phi^*(f) \text{ is finite, and } \phi^*(f)=\int \varphi^*(f)~dP=
\int\left[
f{\varphi'}^{-1}(f)-\varphi\left({\varphi'}^{-1}(f)\right)\right]~dP.
\end{equation}
\end{proposition}

\noindent\textbf{Proof of Proposition  \ref{culcul de phi*}} For
all $f$ in $\langle\mathcal{F}\cup\mathcal{B}_b\rangle$, define
the mapping $G_{f}~:~\mathcal{M}_{\mathcal{F}}~\rightarrow
~]-\infty,+\infty]$~ by
\begin{equation*}
G_{f}(Q):=\phi (Q,P)-\int f~dQ,
\end{equation*}
from which $\phi^*(f)=-\inf_{Q\in
\mathcal{M}_{\mathcal{F}}}G_{f}(Q)$. The function $G_{f}(.)$ is
strictly convex. Its domain is
\begin{equation*}
\text{dom}G_f:=\left\{ Q\in \mathcal{M}_{\mathcal{F}}~\text{ such
that }~G_{f}(Q)<+\infty\right\}.
\end{equation*}
Denote  by  $Q_{0}:=\arg\inf_{Q\in
\mathcal{M}_{\mathcal{F}}}G_{f}(Q)$, which belongs to
$\text{dom}G_f$, if it exists. It follows that $Q_{0}$ is a.c.
w.r.t. $P$. Since $\mathcal{M}_{\mathcal{F}}$ is convex set, the
measure $Q_{0}$ (if it exists) is the only measure in
$\text{dom}G_f$ such that for any measure $R$ in $\text{dom}G_f$,
\begin{equation*}
G_{f}^{\prime}(Q_{0},R-Q_{0})\geq 0,
\end{equation*}
where $G_{f}^{\prime}(Q_{0},R-Q_{0})$ is the directional
derivative of the function $G_{f}$ at point $Q_{0}$ in direction
$R-Q_{0}$; see e.g. Theorem III.31 in \cite{Aze1997}. Denote
$r:=\frac{dR}{dP}$ and $q_{0}:=\frac{dQ_{0}}{dP}$. By its very
definition, we have
\begin{eqnarray*}
G_{f}^{\prime}(Q_{0},R-Q_{0}) & := & \lim_{\epsilon \downarrow
0}\frac{1}{\epsilon}\left\{G_{f}(Q_{0}+\epsilon (R-Q_{0}))-G_{f}(Q_{0})\right\}\\
& = & \lim_{\epsilon \downarrow 0}\int \frac{1}{\epsilon }\left[
\varphi \left( q_{0}+\epsilon (r-q_{0})\right) -\varphi
(q_{0})\right] ~dP-\int f~d(R-Q_{0}).
\end{eqnarray*}
Define the function
\begin{equation*}
g(\epsilon ):=\frac{1}{-\epsilon }\left[ \varphi
\left(q_{0}+\epsilon (r-q_{0})\right) -\varphi (q_{0})\right] .
\end{equation*}
Convexity of $\varphi$ implies
\begin{equation*}
g(\epsilon )\uparrow \varphi^{\prime }(q_{0})(q_{0}-r)~\text{ when
}~\epsilon \downarrow 0,
\end{equation*}
and for all $0<\epsilon\leq 1$ and $R$ in $\text{dom}G_f$, we have
\begin{equation*}
g(\epsilon)\geq g(1)=-\left(\varphi(r)-\varphi(q_{0})\right) \in
L_1(\mathcal{X},P).
\end{equation*}
So, applying  the monotone convergence theorem, we obtain
\begin{equation*}
G_{f}^{\prime}(Q_{0},R-Q_{0})=\int \left(
\varphi^{\prime}(q_{0})-f~\right)~ d(R-Q_{0})\geq 0.
\end{equation*}
Therefore, under assumption (\ref{condition pr calculer phi*}),
for any function $f$  in
$\langle\mathcal{F}\cup\mathcal{B}_b\rangle$ such that
$\text{Im}f\subset \text{Im}\varphi'$, the measure $Q_0$ exists
and it is given by $dQ_0={\varphi'}^{-1}(f)~dP$.  It follows that
\begin{equation}\label{F-L-transform}
\phi^*(f)=\int \varphi^*(f)~dP=\int \left[
f{\varphi'}^{-1}(f)-\varphi\left({\varphi'}^{-1}(f)\right)\right]~dP.
\quad \blacksquare\\
\end{equation}

\begin{remark}\label{remarque concernant int phi*}
If the convex function $Q\in \mathcal{M}_\mathcal{F} \mapsto
\phi(Q,P)\in[0,+\infty]$ is proper, i.e.,
\begin{equation}\label{the function is proper}
\text{there exists at least one measure } Q_0 \text{ in }
\mathcal{M}_\mathcal{F}(P) \text{ such that } \phi(Q_0,P) \text{
is finite},
\end{equation}
then the integral $\int \varphi^*(f)~dP$ is well defined for all
$f\in
 \langle\mathcal{F}\cup\mathcal{B}_b\rangle$. Indeed, for all $f\in
 \langle\mathcal{F}\cup\mathcal{B}_b\rangle$ and for all
 $x\in\mathcal{X}$, by Fenchel's inequality, we have
 \begin{equation}\label{ineg 3}
\varphi^*\left(f(x)\right)\geq
f(x)\frac{dQ_0}{dP}(x)-\varphi\left(\frac{dQ_0}{dP}(x)\right).
 \end{equation}
The RHS term belong to $L_1(\mathcal{X},P)$ by assumption
(\ref{the function is proper}). Hence,  the integral $\int
\varphi^*(f)~dP$ is well defined. Moreover, we have $-\infty <
\int \varphi^*(f)~dP \leq +\infty$ for all
$f\in\langle\mathcal{F}\cup\mathcal{B}_b\rangle$. Hence, from
Theorem \ref{duality for divergences} we can state the following
result.\\
\end{remark}

\begin{theorem}\label{representation duale v2} Assume that $\varphi$ is differentiable.
 Then, for all $Q\in\mathcal{M}_\mathcal{F}$ such that $\phi(Q,P)$ is
 finite and $\varphi'\left(\frac{dQ}{dP}\right)$ belongs to
 $\langle\mathcal{F}\cup\mathcal{B}_b\rangle$, the $\phi$-divergence $\phi(Q,P)$ admits the dual
 representation
\begin{equation}\label{representation duale v2 de phi}
\phi(Q,P)=\sup_{f\in\langle\mathcal{F}\cup\mathcal{B}_b\rangle}\left\{\int
f~dQ - \int \varphi^*(f)~dP\right\},
\end{equation}
and the function
$\overline{f}:=\varphi'\left(\frac{dQ}{dP}\right)$ is a dual
optimal solution. Furthermore, if $\varphi$ is essentially smooth,
then $\overline{f}$ is the unique dual optimal solution
($P$-a.e.).\\
\end{theorem}

\noindent\textbf{Proof of Theorem \ref{representation duale v2}}
 Let $Q\in\mathcal{M}_\mathcal{F}$ such that $\phi(Q,P)$
is finite. Then, the integral $\int \varphi^*(f)~dP$ is well
defined for all $f\in\langle\mathcal{F}\cup\mathcal{B}_b\rangle$;
see Remark \ref{remarque concernant int phi*}. Furthermore, using
(\ref{ineg 3}) for all $Q\in\mathcal{M}_\mathcal{F}(P)$, we can
see that $\phi^*(f)\leq \int \varphi^*(f)~dP$ for all
$f\in\langle\mathcal{F}\cup\mathcal{B}_b\rangle$. Hence, using
Theorem \ref{duality for divergences}, we can write
\begin{equation*}
\phi(Q,P)=\sup_{f\in\langle\mathcal{F}\cup\mathcal{B}_b\rangle}\left\{\int
f~dQ - \phi^*(f)\right\}\geq
\sup_{f\in\langle\mathcal{F}\cup\mathcal{B}_b\rangle}\left\{\int
f~dQ - \int \varphi^*(f)~dP\right\}.
\end{equation*}
On the other hand, by (\ref{eqn 1}), we obtain for the function
$\overline{f}:=\varphi'\left(dQ/dP\right)$,
\begin{equation*}
\varphi^*(\overline{f})=\varphi'\left(\frac{dQ}{dP}\right)\frac{dQ}{dP}-
\varphi\left(\frac{dQ}{dP}\right).
\end{equation*}
From this, using the fact  that the integrals $\int
\left|\varphi'\left(\frac{dQ}{dP}\right)\right|~d|Q|$ and $\int
\varphi\left(\frac{dQ}{dP}\right)~dP$ are finite, by simple
calculus we obtain the equality $\int \overline{f}~dQ - \int
\varphi^*(\overline{f})~dP=\phi(Q,P)$, which completes the proof.
$\blacksquare$\\
\\
\noindent Theorem \ref{representation duale v2} remains valid if
we substitute the vector space
$\langle\mathcal{F}\cup\mathcal{B}_b\rangle$ by the arbitrary
class of function $\mathcal{F}$. We state this result in the
following Theorem.\\

\begin{theorem}\label{representation duale v3} Assume that $\varphi$ is differentiable.
Let $\mathcal{F}$ be an arbitrary class of measurable real valued
functions on $\mathcal{X}$. Then, for all
$Q\in\mathcal{M}_\mathcal{F}$ such that $\phi(Q,P)$ is
 finite and $\varphi'\left(\frac{dQ}{dP}\right)$ belongs to
 $\mathcal{F}$, the $\phi$-divergence $\phi(Q,P)$ admits the dual
 representation
\begin{equation}\label{representation duale v2 de phi}
\phi(Q,P)=\sup_{f\in\mathcal{F}}\left\{\int f~dQ - \int
\varphi^*(f)~dP\right\},
\end{equation}
and the function
$\overline{f}:=\varphi'\left(\frac{dQ}{dP}\right)$ is a dual
optimal solution. Furthermore, if $\varphi$ is essentially smooth,
then $\overline{f}$ is the unique dual optimal solution
($P$-a.e.).\\
\end{theorem}

\begin{remark}
Theorem \ref{representation duale v3}, with an appropriate choice
of the class $\mathcal{F}$, has been used by \cite{Keziou2003} and
\cite{Broniatowski-Keziou2003} to introduce an new common
definition of the ``minimum $\phi$-divergence estimates'' in
discrete or continuous parametric models. Note that the
``plug-in'' minimum $\phi$-divergence estimates introduced by
\cite{Liese-Vajda1987} in chapter 10 are defined only in discrete
parametric models, see also \cite{Lindsay1994} and
\cite{Morales-Pardo-Vajda1995}. The use of the dual representation
(\ref{representation duale v2 de phi}) allows to give a common
definition of the minimum $\phi$-divergence estimates in discrete
or continuous parametric models.\\
\end{remark}

\begin{remark}
Other versions of dual representations of $\phi$-divergences are
given in \cite{BorweinLewis1991} on $L_k(\mathcal{X},P)$ spaces,
in \cite{BorweinLewis1993} on compact metric spaces,  and in
\cite{Leonard2001a} on Orlicz spaces. See also
\cite{Rockafellar1968} for other convex integral functionals on
some ``decomposable'' spaces.\\
\end{remark}

\noindent Under the  assumption
\begin{equation}\label{P integrabilite}
\int_\mathcal{X} |f|~dP \text{ is finite for all }
f\in\mathcal{F},
\end{equation}
the convex function $Q\in
\left[\mathcal{M}_\mathcal{F}(P);\tau_\mathcal{F}\right] \mapsto
\phi(Q,P)\in [0,+\infty]$ is proper. Its Fenchel-Legendre
transform is
\begin{equation}\label{phi* de f v2}
f\in \left[\langle\mathcal{F}\cup\mathcal{B}_b\rangle;
\tau_\mathcal{M}\right] \mapsto
\phi^*(f):=\sup_{Q\in\mathcal{M}_\mathcal{F}(P)}\left\{\int f
~dQ-\phi(Q,P)\right\}\in (-\infty,+\infty],
\end{equation}
which is convex and lower semi-continuous. Following
\cite{Rockafellar1968} p. 532, let
$L^*:=\langle\mathcal{F}\cup\mathcal{B}_b\rangle$ and
$L:=\mathcal{M}_\mathcal{F}(P)$. Then condition (\ref{P
integrabilite}) implies that both $L^*$ and $L$ are decomposable.
 Hence, we can apply the
Corollary of Theorem 2 in \cite{Rockafellar1968}, to obtain the
following result:\\

\begin{theorem}\label{representation duale}
Under assumption (\ref{P integrabilite}), the convex conjugate
function $f\in \langle\mathcal{F}\cup\mathcal{B}_b\rangle \mapsto
\phi^*(f)$  defined by (\ref{phi* de f v2}) is proper, it can be
expressed by
\begin{equation}
  \phi^*(f)=\int \varphi^*(f)~dP \text{ for all } f\in
  \langle\mathcal{F}\cup\mathcal{B}_b\rangle,
\end{equation}
and the $\phi$-divergence $\phi(Q,P)$ admits the dual
representation
\begin{equation}
 \phi(Q,P) = \sup_{f\in
 \langle\mathcal{F}\cup\mathcal{B}_b\rangle}\left\{\int f~dQ - \int
 \varphi^*(f)~dP\right\}, \text{ for all } Q\in
 \mathcal{M}_\mathcal{F}(P).
\end{equation}
In particular, the function
$Q\in\left[\mathcal{M}_\mathcal{F}(P);\tau_\mathcal{F}\right] \to
\phi(Q,P)\in [0,+\infty]$ is lower semi-continuous.\\
\end{theorem}

\begin{remark}
The lower semi-continuity property of the function
$$Q\in\left[\mathcal{M}_\mathcal{F}(P);\tau_\mathcal{F}\right]
\mapsto \phi(Q,P)\in [0,+\infty]$$ holds from Proposition
\ref{prop div lsc param} and Lemma \ref{Lemma M-P-F ferme} without
assuming (\ref{P integrabilite}). On the other hand, Theorem
\ref{representation duale v3} and \ref{representation duale} are
of interest particularly when $\phi(Q,P)$ is finite  and the class
$\mathcal{F}$ contains the function $f=\varphi'(dQ/dP)$. In
Theorem \ref{representation duale v3}, condition on $Q$, i.e.,
$\int |\varphi'(dQ/dP)|~ d|Q| <\infty$,  holds whenever
$\phi(Q,P)$ is finite and $\varphi$ satisfies condition
(\ref{condition C.0}); see Lemma \ref{lemma caract proj}. However,
in Theorem \ref{representation duale}, these conditions do not
inevitably imply assumption (\ref{P integrabilite}) if the class
$\mathcal{F}$ contains $\varphi'(dQ/dP)$. It is the case, for
example, when $\phi=KL$, $Q$ is a normal law and $P$ is a Cauchy
law. Indeed, $KL(Q,P)$ is finite, the assumption $\int
|\log(dQ/dP)|~dQ < \infty$ in Theorem \ref{representation duale
v3} holds  while the assumption  $\int |\log(dQ/dP)|~dP < \infty$
in Theorem \ref{representation duale} does not. This  shows the
interest of Proposition \ref{prop div lsc param} and Theorem
\ref{representation duale v3}.
\end{remark}

\section{Applications to the minimization of $\phi$-divergences
on sets of signed finite measures satisfying linear constraints}
 \noindent In this section we apply the results of the sections 2, 3 and 4 to
the optimization problem
\begin{equation*}
\inf_{Q\in M_g}\phi(Q,P)
\end{equation*}
where $M_g$ is defined in (\ref{M sets}).\\
\\
\noindent Under different assumptions, we obtain the dual equality
$\inf(\ref{optim prob})=\sup(\ref{lagrangian dual prob})$ and
results about the problems of existence, uniqueness and
characterization of the dual optimal solution and
the $\phi$-projections  of $P$ on the set $M_g$. \\
\\
\noindent We state our results under the following assumptions:
\begin{eqnarray}
 & & \text{the convex function } \varphi \text{ is differentiable};\label{condition de
 differ} \\
 & &  \text{there exists at least one } \phi\text{-projection } Q^*
\text{ of }  P \text{ on } M_g   \text{ with the same support as }
P. \label{condition pour ldp 1 M}
\end{eqnarray}
\begin{theorem}\label{ecriture duale M}
Assume that conditions (\ref{condition C.0}), (\ref{condition de
differ}) and (\ref{condition pour ldp 1 M})
 hold. Then
\begin{enumerate}
\item [(1)] there exists $\overline{\lambda}\in \mathbb{R}^{1+l}$
such that
\begin{equation*}
\varphi'\left(\frac{dQ^*}{dP}(x)\right)=\overline{\lambda}_0+\sum_{i=1}^{l}\overline{\lambda}_ig_i(x)\quad
(P-a.e.),
\end{equation*}
\item [(2)] the  equality
\begin{equation*}\label{ecriture duale 2 M}
\inf_{Q\in M_g}
\phi(Q,P)=\sup_{\lambda\in\mathbb{R}^{1+l}}\left\{\lambda_0-\int_\mathcal{X}
\varphi^*\left(\lambda^Tg(x)\right)~dP(x)\right\}
\end{equation*}
holds, and $\overline{\lambda}$ is a dual optimal solution.
Furthermore, if the function $\varphi$ is essentially smooth, then
the dual optimal solution $\overline{\lambda}$ is unique.\\
\end{enumerate}
\end{theorem}

\begin{remark} Under assumptions of Theorem \ref{ecriture duale
M}, the $\phi$-projection of $P$ on $M_g$ is characterized without
supposing that $\overline{\lambda}$ is an interior point of
$\text{dom}\phi^*$. Furthermore, the dual equality holds and the
dual optimal solution is attained. Sufficient conditions for
assumption (\ref{condition pour ldp 1 M}) are given in Corollary
\ref{corollaire 1} and   Proposition \ref{proposition 3} below.\\
\end{remark}

\noindent \textbf{Proof of Theorem \ref{ecriture duale M}} Under
assumptions (\ref{condition C.0}) and (\ref{condition de differ}),
part (1) is a direct consequence of Theorem \ref{thm rusch} part
(2). We prove now part (2). We have $\inf_{Q\in M_g}
\phi(Q,P)=\phi(Q^*,P)$ since $Q^*$ is a $\phi$-projection of $P$
on $M_g$. Now, by Theorem \ref{representation duale v3}, choosing
the class of measurable functions
$$\mathcal{F}=\left\{x\in\mathcal{X}\mapsto \lambda^Tg(x) \text{
such that } \lambda\in\mathbb{R}^{1+l}\right\},$$ we can write
\begin{equation*}
 \phi(Q^*,P)=\sup_{\lambda\in\mathbb{R}^{1+l}}
 \left\{\lambda_0-\int_\mathcal{X}\varphi^*\left(\lambda^Tg(x)\right)~dP(x)\right\},
\end{equation*}
and from it we deduct that $\overline{\lambda}$  is a dual optimal
solution by the same Theorem. $\blacksquare$\\

\begin{corollary}\label{corollaire 1}
Assume that $\varphi$ is differentiable and strictly convex. If
there exists some $\overline{\lambda}\in\mathbb{R}^{1+l}$ such
that
\begin{equation}\label{condition reciproque}
 \int
\varphi\left({\varphi'}^{-1}\left({\overline{\lambda}}^Tg(x)\right)\right)~dP<\infty
\quad \text{and} \quad\int{g}^T
{\varphi'}^{-1}\left({\overline{\lambda}}^Tg(x)\right)~dP(x)=(1,0,\ldots,0)^T,
\end{equation}
then
\begin{enumerate}
 \item [(1)] the measure $Q^*$ defined by
$dQ^*(x)={\varphi'}^{-1}\left({\overline{\lambda}}^Tg(x)\right)~dP(x)$
is the unique $\phi$-projection of $P$ on $M_g$.
 \item [(2)] the equality
\begin{equation*}\label{ecriture duale 2 M}
\inf_{Q\in M_g}
\phi(Q,P)=\sup_{\lambda\in\mathbb{R}^{1+l}}\left\{\lambda_0-\int_\mathcal{X}
\varphi^*\left(\lambda^Tg(x)\right)~dP(x)\right\}
\end{equation*}
holds, and $\overline{\lambda}$ is a dual optimal solution.
Furthermore, if the function $\varphi$ is essentially smooth, then
the dual optimal solution $\overline{\lambda}$ is unique.
\end{enumerate}
In particular,  (\ref{condition reciproque}) holds if there exists
a dual optimal solution $\overline{\lambda}$ which is an interior
point of
\begin{equation*}
\text{dom} \phi^* :=\left\{\lambda\in\mathbb{R}^{1+l}~\text{ such
that }\int_\mathcal{X}\left|
\varphi^*\left(\lambda^Tg(x)\right)\right|~dP(x) \text{
is finite }\right\}.\\
\end{equation*}
\end{corollary}

\noindent \textbf{Proof of Corollary \ref{corollaire 1}} (1) Apply
Theorem \ref{thm rusch} part (1). (2) the proof is the same as
that of part (2) of Theorem \ref{ecriture
duale M}. $\blacksquare$\\
\begin{remark}
Note that, if $\varphi$ is differentiable and strictly convex,
then under  assumption (\ref{condition C.0}),  conditions
(\ref{condition pour ldp 1 M}) and (\ref{condition reciproque})
are equivalent; see Theorem \ref{thm rusch} part (1) and (2).\\
\end{remark}

\noindent In the the following Proposition we give other
sufficient conditions for assumption (\ref{condition pour ldp 1
M}). The conditions are
\begin{eqnarray}
& & \phi\left(M_g,P\right)<\infty;\label{C.1}\\
& & \lim_{|x|\to\infty}\frac{\varphi(x)}{|x|}=+\infty;\label{C.2}\\
& & \text{for every } \alpha>0, \text{ and  all } i=1,\ldots,l,~
\int \varphi^*\left(\alpha|g_i|\right)~dP<\infty;\label{C.3}\\
& & \text{there exists numbers } 1< r, k < +\infty \text{ such
that } r^{-1}+k^{-1}=1,\label{C.3.1} \\
& &  \lim_{\left\vert x\right\vert \rightarrow \infty
}\frac{\varphi(x)}{\left\vert x\right\vert ^{r}}>0, ~\text{ and
for all } i=1,\ldots,l, ~ \left\|g_i\right\|_k<\infty;\nonumber\\
& & \text{the functions } g_1,\ldots,g_l  \text{ belong to  }
L_\infty(\mathcal{X},P);\label{C.4}\\
& & \varphi(0)=+\infty;\label{C.5}\\
& & a_\varphi=0 \quad \text{and} \quad \varphi'(0)=-\infty;
\label{C.6} \\
& & \text{there exits some }  Q_0\in M\cap \text{dom}\phi \text{
such that } \frac{dQ_0}{dP} > 0 ~~ (P-a.e.).\label{C.7}
\end{eqnarray}

\begin{proposition}\label{proposition 3}\
\begin{enumerate}
 \item [(1)] Under assumptions (\ref{C.1}), (\ref{C.2}), (\ref{C.3}) and
(\ref{C.5}), condition (\ref{condition pour ldp 1 M}) holds.
 \item [(2)] Condition (\ref{condition pour ldp 1 M}) holds also under assumptions
(\ref{C.1}), (\ref{C.3.1}) and (\ref{C.5}).
  \item  [(3)] Condition
(\ref{condition pour ldp 1 M}) holds also if, in part (1),
(\ref{C.3}) is replaced by (\ref{C.4}) or/and if condition
(\ref{C.5}) is replaced by $\left[(\ref{condition de differ}),
(\ref{C.6}) \text{ and } (\ref{C.7})\right]$.
 \item [(4)] Condition
(\ref{condition pour ldp 1 M}) holds also if, in part (2),
condition (\ref{C.5}) is replaced by $\left[(\ref{condition de
differ}), (\ref{C.6})\right.$  and  $\left.(\ref{C.7})\right]$.\\
\end{enumerate}
\end{proposition}

\noindent \textbf{Proof of Proposition \ref{proposition 3}} (1)
Since $M_g$ is closed in
$\left[\mathcal{M}_\mathcal{F};\tau_\mathcal{F}\right]$ (choosing
the class $\mathcal{F}=\left\{g_1,\ldots,g_l\right\}$), we can
then apply Theorem \ref{thm exist proj sur ferme 1} to deduce that
there exists at least one $\phi$-projection of $P$ on $M_g$.
Condition (\ref{C.5}) implies that $Q^*$ has the same support as
$P$. (2) We can apply Theorem \ref{thm exist proj sur ferme}. (3)
Under assumption (\ref{C.4}), the set $M_g$ is closed in
$\tau$-topology. Hence, we can apply Theorem \ref{corollary thm
exist proj sur ferme} to deduce that there exists at least one
$\phi$-projection of $P$ on $M_g$. Conditions (\ref{condition de
differ}), (\ref{C.6}) \text{ and } (\ref{C.7}) imply that $Q^*$
has the same support as $P$ (see Lemma \ref{meme
support}). $\blacksquare$\\

\bibliographystyle{natbib}

\end{document}